\documentclass[11pt,a4paper]{article}

\usepackage[T1]{fontenc}
\usepackage[utf8]{inputenc}
\usepackage{mathptmx}
\usepackage{helvet}
\usepackage{courier}

\usepackage{graphicx}        % standard LaTeX graphics tool
\usepackage{multicol}        % used for the two-column index
\usepackage[bottom]{footmisc}% places footnotes at page bottom
\usepackage{amssymb, amsmath, latexsym}%amsthm,
\usepackage{mathrsfs}
\usepackage{eurosym}
\usepackage{paralist}
\usepackage[margin=1in]{geometry}
\usepackage{amsfonts}
\usepackage{amssymb,enumerate}
\newcommand{\qed}{}
\usepackage{breqn}

\usepackage{relsize}

\allowdisplaybreaks

\makeindex             % used for the subject index
\newenvironment{proof}[1][Proof]%
{\vspace{0.5em}\par\textit{{#1}.}\;\;}{\phantom{|}\hfill$\Box$\par}

\begin{document}

\title{Subdominant eigenvalue location and
the\\
 robustness of Dividend Policy Irrelevance}

\author{Adam J. Ostaszewski}
\maketitle

\textbf{Abstract. } This paper, on subdominant eigenvalue location of a bordered diagonal matrix, is the mathematical sequel to an accounting
paper by Gao, Ohlson, Ostaszewski \cite{GaoOO}. We explore the following
characterization of dividend-policy irrelevance (DPI) to equity valuation in
a multi-dimensional linear dynamics framework $L$: DPI occurs under $L$ when
discounting the expected dividend stream by a constant interest rate iff
that rate is equal to the dominant eigenvalue of the canonical principal
submatrix $A$ of $L.$ This is justifiably the `latent' (or gross) rate of return, since
the principal submatrix relates the state variables to each other but with
dividend retention. We find that DPI reduces to the placement of the
maximum eigenvalue of $L$ between the dominant and subdominant eigenvalues
of $A.$ We identify a special role, and a lower bound, for the coefficient measuring the
year-on-year dividend-on-dividend sensitivity in achieving robust equity
valuation (independence of small variations in the dividend policy).

\section{Introduction and motivation}

Accounting theory seeks to reconcile valuation of a firm based on \textit{%
historically} observed variables (`primitives', that recognize value created
\textit{to date}) with its equity value, arrived at by markets in a \textit{%
prospective } fashion. The market's valuation is theoretically modelled as
the present value of future (expected) dividends and involves discounting by
the (notional) \textit{riskless} interest rate in force, say $r$ per unit
time. From the historic (accounting) side, various secondary composite
variables have been derived from the primitives (with appropriate technical
names such as `residual income' -- for a brief introduction see \cite{GieO}), formalizing in one way or
another a notion of \textit{current} `earnings'; the latter is then
intended to identify equity value directly (as a dependent variable) and to
provide empirically stable time series.

To arrive at such a composite accounting variable, assumptions are needed
concerning the future evolution of the primitives -- at least in a
hypothetical `steady state' context. (For a `dynamic' alternative, drawing on the value of waiting, see \cite{DavO} in this same volume.) The favourite mechanism for this context is a
\textit{linear state-space representation }$L$, thereby introducing subtle
links -- our main concern here -- between accounting theory and mathematics.

It is noteworthy, though not of direct mathematical significance to this
paper, that an encouraging feature for the use of a (linear) representation $%
L$ is its flexibility in permitting inclusion, alongside state variables
that recognize historic value creation (as above), additional `information'
state variables; these capture the (typical) dynamics of an embedded
`potential to create' value, an `intangible' value, currently unrecognized
in the accounts but feeding through to future recognized value (a matter
central to the luckless 2014 attempt by Pfizer to bid for AstraZeneca --
`the mega-merger that never was'). This partly bridges the
historic-prospective divide. (The idea was introduced into the accounting
literature of linear systems by Ohlson \cite{Ohl95}, and enabled him to include
the accounting of `goodwill' value -- see \cite{LoL}; for another example of an
intangible, involving product `image' and its valuation, see e.g. \cite{JacJZ}.)

Returning to mathematical concerns, we note that the eigenstructure of $L$
(eigenvalue distribution) has to connect with economic consequences of an
assumed `steady state' -- such as absence of \textit{arbitrage}
opportunities in equity valuation, and its relation to the notional riskless
interest rate (above). A further fundamental insight, going back to Miller
and Modigliani \cite{MilM2} in 1961, is that -- under prescribed conditions (but
see e.g. \cite{Bha} for the effects of alternative informational assumptions) --
the equity value should not depend on the distribution of value, be it
impounded into the share price or placed in the share-holders' pockets (via
dividend payouts); this is properly formalized below. (This is one aspect of
capital structure irrelevance: equity value should not depend on debt versus
equity issuance \cite{MilM1,MilM3}.) The principle of \textit{dividend policy
irrelevancy} also carries implications for the linear dynamics. For a recent
analysis of the connections see \cite{GaoOO}, where the basic result asserts that
DPI\ occurs iff the riskless interest rate agrees with the dominant
eigenvalue of the reduced linear system (`subsystem') obtained by the firm
withholding (retaining) dividend payouts. One may call the latter the
dominant `latent' rate of the system $L$. Recall that it is the riskless
rate $r$ that is used in the present-value calculation above.

This is a \textit{knife-edge} characterization in regard both to the
riskless interest rate and the dominant eigenvalue, so it is natural to
study \textit{accounting robustness} in the DPI\ framework. That is the
principal aim of this paper, achieved by studying an eigenvalue location
problem, similar but distinct from one in control theory (reviewed shortly
below). The delicacy of this matter is best seen in the light of Wilkinson's
example in \cite[\S 33]{Wil} of a sparse $20\times 20$ matrix with the integers $%
n=1,2,...,20$ on the diagonal (its eigenvalues) and all superdiagonal
entries of $20$; a small perturbation of $\varepsilon $ in the bottom
left-hand corner yields the characteristic polynomial to be $(\lambda
-20)...(\lambda -1)-20^{19}\varepsilon $, so that for $\varepsilon
:=10^{-10} $ the eigenvalues are these: 6 real ones which are to 1 decimal
place 0.9, 2.1, 2.6, 18.4, 18.9, 20.0, and 7 conjugate complex pairs (all in
mid-range values, in modulus) as follows:

\begin{equation*}
4\pm i1.1;\;5.9\pm i1.9;\;8.1\pm i2.5;\;10.5\pm i2.7;\;12.9\pm i2.5;\;15.1\pm
i1.9;\;17\pm i1.1.
\end{equation*}

Turning now to the mathematical problem, consider, granted initial
conditions, the performance of the following discrete-time system:
\begin{equation}
\left.
\begin{array}{cccc}
z_{t+1}\ = & Az_{t} & +\ av_{t}\  & +bd_{t}, \\
v_{t+1}\ = &  & +\ \alpha v_{t}, &  \\
d_{t+1}\ = & w^{T}z_{t} & +\ \gamma v_{t} & +\ \beta d_{t}.%
\end{array}%
\right\}  \tag{$\Omega $}
\end{equation}%
Here $\beta $ is the \textit{dividend-on-dividend} `year-on-year' growth;
its effect is particularly significant -- see below and \S 2 (Theorem 3). So
the state variables at time $t$ are $(z_{t},v_{t},d_{t})\in \mathbb{R}^{n+2}$
-- with $d_{t}$ representing the time-$t$ dividend, $v_{t}$ an `information'
variable (as above), subjected to `fading' over time by a factor $\alpha $
satisfying%
\begin{equation*}
0\leq \alpha <1,
\end{equation*}%
with $A$ a real matrix (hereafter, the \textit{reduced matrix} of the
system, or the `dividend-retention' matrix) that is constant over time. The
performance of the system at time $0$ is measured by the expression%
\begin{equation*}
P_{0}=\sum_{t=1}^{\infty }R^{-t}d_{t},
\end{equation*}%
in which a discount factor $R$ is applied to the sequence $d=\{d_{t}\}$
generated by $(\Omega ).$ Here $R=1+r$ with $r$ as above (the governing
riskless interest rate per unit time), and so $P_{0}$ represents an initial
\textit{equity} valuation of the firm -- in the sense motivated above. To
guarantee convergence it is sufficient to assume that all components of any
solution of $(\Omega )$ have growth below $R;$ referring to the modulus of
the dominant eigenvalue of the coefficient matrix in $(\Omega )$ by $\lambda
_{\text{max}}^{\Omega }$ this growth condition may be restated as%
\begin{equation}
\lambda _{\text{max}}^{\Omega }<R.  \label{lambda-max}
\end{equation}%
The bottom row vector in $(\Omega )$%
\begin{equation*}
\omega _{\text{div}}=(w^{T},\gamma ,\beta )
\end{equation*}%
is termed the \textit{dividend policy}.

In this setting the full coefficient matrix is assumed constant, but not
known to observers of the state variables (which are disclosed in the annual accounts).
However, whereas $A$ and the penultimate row involve value created over time
through an initial (fixed) investment, it is the final row that generates
the returns (over time) to the investors. Thus the equity $P_{0}$ should be
regarded as a function of $R$ and of the dividend-policy vector parameters
set by the managers, that is%
\begin{equation}
P_{0}(R;\omega _{1},...,\omega _{n},\gamma, \beta):=\sum_{t=1}^{\infty
}R^{-t}d_{t}.  \label{equity}
\end{equation}%
One says that the valuation $P_{0}(.)$ exhibits \textit{Dividend Policy
Irrelevancy} (DPI)\textit{\ at} $R$ if the function $P_{0}(R;...)$ is
unchanged as the dividend-policy vector $\omega _{\text{div}}$ varies. A
first problem is to determine circumstances under which the system exhibits
dividend irrelevancy. Up to a technical side-condition (ensuring
co-dependence of dividends and value creation) the short answer is that the
\textit{dominant eigenvalue} of $A$ should agree with $R$ -- this was first
proved by Ohlson in the special case $n=1,$ and then generalized in \cite{GaoOO}
(and also referred to in the earlier published monograph \cite{OhlG}).

Below we refine the notion of dividend irrelevancy in order to study the
effects of a proximal \textit{sub-dominant} eigenvalue. We first establish
notation and some conventions. Begin by omitting hereafter explicit mention
of the information variable $v_{t};$ we regard it as yet another state
variable absorbed into $A$ (with $\alpha $ then becoming an eigenvalue of
the reduced matrix), and so we overlook its simple dynamics; we may now free
up $\alpha $ and $\gamma$ for other uses below. The eigenvalues $\lambda _{1}^{A},\lambda
_{2}^{A},...,\lambda _{n}^{A}$ of the reduced matrix $A$ (\textit{latent}
relative to $L$, below) are listed in order of decreasing modulus. As these
will be required to be real, positive and (generically) distinct, this is
taken to mean%
\begin{equation*}
\lambda _{1}^{A}>\lambda _{2}^{A}>...>\lambda _{n}^{A}>0.
\end{equation*}%
Whenever convenient (e.g. in proofs) we omit the superscript $A.$ The system
matrix $L$ of ($\Omega $) above, viewed as the \textit{augmented matrix} of $%
A$, is now given by%
\begin{equation*}
L_{A}=A^{\lrcorner }(w,\beta ):=\left[
\begin{array}{cc}
A & b \\
w^{T} & \beta%
\end{array}%
\right]
\end{equation*}%
("$A$-bordered"), and is regarded as a function of the \textit{real} vector $(w^{T},\beta ).$
Its (possibly complex) eigenvalues will likewise be regarded as \textit{%
functions} of $(w^{T},\beta )$ and denoted by $\lambda _{j}^{L},$ or more
simply by $\kappa _{j},$ so that
\begin{equation*}
\kappa _{1}=\lambda _{1}^{L},\kappa _{2}=\lambda _{2}^{L},...,\kappa
_{n+1}=\lambda _{n+1}^{L},
\end{equation*}%
which distinguish them more easily (from $\lambda _{j}^{A});$ each index
here is identified through the functional conditions
\begin{equation}
\kappa _{j}(0,...,0,\beta )=\lambda _{j}^{A}\text{ for }j=1,...,n,\text{ and
}\kappa _{n+1}(0,...,0,\beta )=\beta .  \label{initialize}
\end{equation}%
We write%
\begin{equation*}
\lambda _{\max }^{L}(w^{T},\beta),\text{ or }\kappa _{\max }(w^{T},\beta
):=\max \{\kappa _{j}(w^{T},\beta ):j=1,...,n+1\}.
\end{equation*}%
Although $A^{\lrcorner }$ is not in general symmetric, we will contrive
situations in which the eigenvalues of $A$ interlace with those of $%
A^{\lrcorner }:$ $\kappa _{j}\geq \lambda _{j}\geq \kappa _{j+1},$ just as
in Cauchy's Interlace Theorem, cf. \cite[Th. 4.3.17]{HorJ}, \cite[Ch.7, \S 8 Th. 4]{Bel}, \cite{Hwa}, at least for $j\geq2$.

Since we are mostly concerned with the characteristic polynomial and
eigenvalue location, we will be working in an equivalent \textit{canonical
setting} in which, firstly, $A$ is diagonal and, secondly, as a further
simplification, we suppose that\textit{\ }for $j\leq n$ the\textit{\
dividend significance coefficients }$b_{j}$\textit{\ are all non-zero. }%
Rescaling by $b_{j}$ the $j$-th equation of the diagonalized system gives
what we term the equivalent \textit{canonical} system in which the resulting
\textit{dividend significance coefficients are }$\delta _{j}=\pm 1$ (as, of
course, we may also rescale by $-b_{j}$: see Remark 2 in \S 2.1). Thus $L$
is replaced by%
\begin{equation}
H(\omega )=\left[
\begin{array}{ccccc}
\lambda _{1}^{A} & 0 &  & 0 & \delta _{1} \\
0 & \lambda _{2}^{A} &  &  & \delta _{2} \\
&  & ... &  &  \\
0 & 0 &  & \lambda _{n}^{A} & \delta _{n} \\
\omega _{1} & \omega _{2} &  &  & \omega _{n+1}%
\end{array}%
\right] ,  \label{canon}
\end{equation}%
where $\delta _{j}=\pm 1$ for each $j.$
It is preferable to subsume $\beta$ as $\omega_{n+1}$ (rather than as $\delta_{n+1}$) into the `canonical dividend-policy vector' corresponding to $(w^T,\beta)$.
Our first definitions all contain
growth conditions analogous to (\ref{lambda-max}) and are motivated by
Proposition 1 below.

\bigskip

\textbf{Definition 1.} \textit{We say that the system }$(\Omega )$\textit{\
has dividend irrelevance at }$R$\textit{\ if, for all }$\omega =(\omega
_{1},...,\omega _{n},\beta )$\textit{\ such that }$|\lambda _{\max
}^{L}(\omega ,\beta )|<R,$%
\begin{equation*}
P_{0}(R;\omega ,\beta )=P_{0}(R;o,\beta ),
\end{equation*}%
\textit{with }$o$\textit{\ the zero vector.}\bigskip

\textbf{Definition 2.} \textit{We say that the system }$(\Omega )$ \textit{%
has local dividend irrelevance at }$R$\textit{\ for }$\omega =(\omega
_{1},...,\omega _{n},\beta )$\textit{\ if there is }$\varepsilon >0$\textit{%
\ so that, for all }$\omega ^{\prime }=(\omega _{1}^{\prime },...,\omega
_{n}^{\prime },\beta ^{\prime })$\textit{\ such that }$||\omega -\omega
^{\prime }||<\varepsilon $\textit{\ and }$|\lambda _{\max }^{L}(\omega
^{\prime },\beta ^{\prime })|<R,$%
\begin{equation*}
P_{0}(R;\omega ^{\prime },\beta ^{\prime })=P_{0}(R;\omega ,\beta ).
\end{equation*}

Here the norm is Euclidean. The local definition is weaker in that it
requires merely that the equity valuation be robust in respect of the
accounting system (i.e. insensitive to minor accounting variations).
However, in our model setting $P_{0},$ regarded as a function of $\omega
_{1},...,\omega _{n+1},$ is a rational function in these variables (see
Observation below in \S 1.2), so its local constancy for a given $R$ is
equivalent to global constancy for the same $R$. An intermediate definition
permitting constant equity is the following

\bigskip

\textbf{Definition 3.} \textit{We say that the system (}$\Omega )$\textit{\
has bounded dividend irrelevance at }$R$\textit{\ if, for some positive }$%
\rho <R$\textit{\ and all }$\omega =(\omega _{1},...,\omega _{n})$\textit{\
such that }$|\lambda _{\max }^{L}(\omega ,\beta )|<\rho ,$%
\begin{equation*}
P_{0}(R;\omega ,\beta )=P_{0}(R;o,\beta ).
\end{equation*}

The example below identifies anomalous behaviour which these definitions
offer as possible.

The requirement for dividend irrelevance amounts to discovering to what
extent $P_{0}(R,d)$ depends only on the initial data: $A,b,z_{0},d_{0}.$

In view of the role of the interest rate $R>0,$ it will be appropriate to
make the following.

\bigskip

\textbf{Blanket assumption. }\textit{The eigenvalues }$\lambda
_{1}^{A},...,\lambda _{n}^{A}$\textit{\ of }$A$\textit{\ are all real and
positive.}

\bigskip

Notice that if $|\kappa _{\max }(\omega ,\beta )|<\lambda _{2},$ small
enough variations in the dividend-policy vector will ensure that the
inequality is preserved. This entails (see Proposition 1 below) that the
system will have local dividend-policy irrelevance at more than one rate,
namely at $R=\lambda _{1}$ and $R=\lambda _{2}.$ Our contribution is to
identify in Theorem 3 below a condition on $\beta ,$ namely that%
\begin{equation*}
\beta >2\lambda _{2}-\lambda _{1},
\end{equation*}%
requiring a lower bound on the dividend-on-dividend yearly growth, which
ensures that $|\kappa _{\max }(\omega ,\beta )|>\lambda _{2}$ and thereby
achieves uniqueness of the latent rate of return in this case: dividend-policy irrelevance occurs only at the one rate
$R=\lambda _{1}$.

\bigskip

\textbf{Example of Bounded-DPI at both }$R=\lambda _{1}$ \textbf{and} $%
\lambda _{2}$\textbf{. }We take $\lambda _{1}=2,\lambda _{2}=1.5$, $\beta
=0.5,$ and $\delta _{1}=-1,\delta _{2}=+1,$ $\omega _{2}=0.1.$ Note that $%
2\lambda _{2}-\lambda _{1}=3-2=1,$ so $\beta <2\lambda _{2}-\lambda _{1}.$

\begin{figure}[!t] %F1
\centering\includegraphics[width=0.7\textwidth] {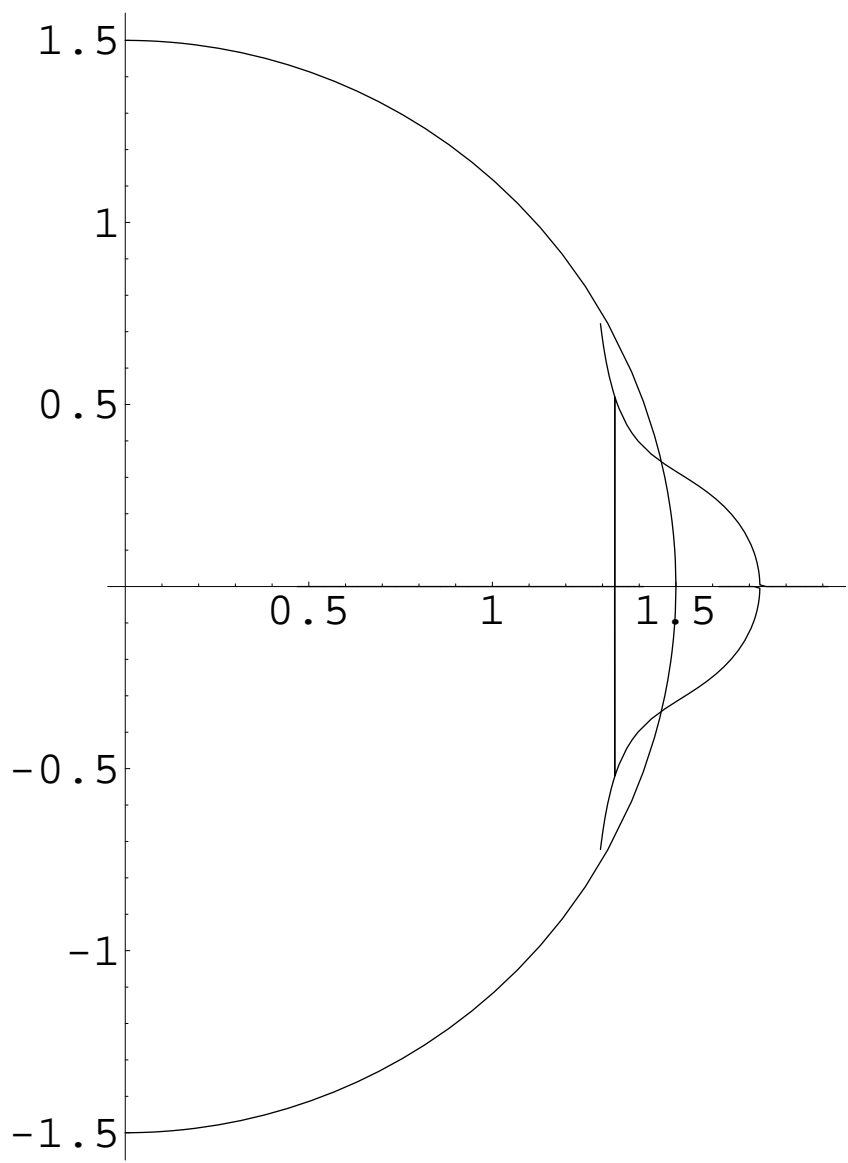}
\caption{
Eigenvalue locus as $\omega_1$ varies. For a range of values of $\omega_1$ all eigenvalues of $A^{\lrcorner }$ lie in the disc of radius $\lambda^A_2$ (the second largest eigenvalue of $A$).
} \label{gig1}
\end{figure}

Figure $\ref{gig1}$ shows the locus of the conjugate complex root pair for $\omega _{1}$
running through the range 0.1 to 1.2 generated by Mathematica. The third
root is below, but close to, $1.5.$ Note that the additional vertical line
appears from the numerics (the computer routine switches the identity of the
two conjugates). Variations in $\omega $ bounding the eigenvalues to $|\zeta
|<\lambda _{2}$ keep both equity $P_{0}(\lambda _{2})$ and $P_{0}(\lambda
_{1})$ constant.

\bigskip

\subsection{Continuous-time analogue}

Letting $D_{t}$ denote cumulative dividends, the analogous dynamic in
continuous time may be formulated as%
\begin{eqnarray*}
\dot{z}_{t} &=&Az_{t}+bD_{t}, \\
\dot{D}_{t} &=&w^{T}z_{t}+w_{n+1}D_{t}.
\end{eqnarray*}%
Here $z_{t}=(x_{t}^{1},...,x_{t}^{n})$ is now the state vector of flow
variables and $b$ and $w$ are columns in $R^{n}$. Putting $%
Z_{t}=(x_{t}^{1},...,x_{t}^{n},D_{t}),$ we arrive at the system%
\begin{equation*}
\dot{Z}_{t}=A^{\lrcorner }Z_{t},
\end{equation*}%
where as before $A^{\lrcorner }$ is $A$ bordered by $b$ and $\bar{w}%
^{T}=(w_{1},...,w_{n},w_{n+1}).$ Thus $D_{t}$ has growth rate at most $%
\kappa _{\max },$ the largest eigenvalue of $A^{\lrcorner }.$

The equity is%
\begin{equation*}
P_{0}=\int_{0}^{\infty }e^{-rt}dD_{t}=-rD_{0}+r\int_{0}^{\infty
}e^{-rt}D_{t}dt,
\end{equation*}%
which converges provided that the dividend flow is of exponential growth
less than $r$ and in particular that%
\begin{equation}
\lim_{t\rightarrow \infty }e^{-rt}D_{t}=0.  \label{nobubble}
\end{equation}%
This requires therefore that $\kappa _{\max }<r.$

The identical format allows interpretation of results in this paper leading
to valid conclusions for the continuous-time framework. In fact set%
\begin{equation*}
\tilde{Z}=\int_{0}^{\infty }e^{-rt}Z_{t}dt
\end{equation*}%
to obtain the discounted flow, i.e. the Laplace transform. We conclude that%
\begin{equation*}
-rZ_{0}+r\tilde{Z}=A^{\lrcorner }\tilde{Z},
\end{equation*}%
subject to $\kappa _{\max }<r.$ Assuming $rI-A^{\lrcorner }$ is invertible
we have%
\begin{equation*}
\tilde{Z}=(rI-A^{\lrcorner })^{-1}Z_{0},
\end{equation*}%
whence%
\begin{equation*}
P_{0}=P_{0}(\bar{w})=(0,...,0,1)^{T}(rI-A^{\lrcorner }(\bar{w}))^{-1}Z_{0}.
\end{equation*}%
This formula was exploited in Ashton et al. in \cite{Ash}, \cite{AshCTW} in a
stochastic setting.

An identical formula can be developed for discrete time (with Laplace
transform replaced by z-transform) and is the departure point for the purely
algebraic argument developed in \cite{GaoOO}. We note an important conclusion,
which follows from the form of the adjugate matrix (or by invoking Cramer's
Rule).\medskip

\textbf{Observation.} In the framework above, the equity $P_{0}$ regarded as
a function of $w_{1},...,w_{n+1}$ is a rational function in these variables.

\subsection{Control-theory analogies}

\ We remark on the contextual similarities between the accounting model and
two standard control-theory settings and identify the differences. For
background, see e.g. \cite{Rug}.

The first setting is the closest. Here we regard $z_{t}$ as a state vector
and $d_{t}$ as a control variable, and write the system equation as%
\begin{equation*}
z_{t+1}=Az_{t}+d_{t}b,
\end{equation*}%
with the control variable selected according to a non-standard (since it is
in effect differential) feedback-law%
\begin{equation*}
d_{t+1}=w^{T}z_{t}+\beta d_{t}.
\end{equation*}%
In these circumstances we interpret $P_{0}(R,d)$ defined by (\ref{equity})
as the system performance index (a discounted cost stream), and require that
that it be either (i) independent of the feedback law parameters, or (ii)
independent of small variations in the parameters. This is akin to
performance stability, or even a `disturbance decoupling problem', except
that there is no modelled disturbance here. Our analysis is similar to that
of the standard pole-placement problem (requiring assignment of eigenvalues
through feedback-law design). Where we differ is in the inevitable presence
and effect of zeros as well as poles.

The alternative view is to regard $z_{t}$ as the observation vector of the
state $Z_{\,t}=(z_{\,t},d_{t}),$ with state equation and observation vector
defined respectively by
\begin{eqnarray*}
Z_{\,t+1} &=&HZ_{\,t}, \\
z_{\,t} &=&PZ_{\,t}.
\end{eqnarray*}%
Here $P$ is the projection matrix from $\mathbb{R}^{n+1}$ to $\mathbb{R}%
^{n}. $ In these circumstances the performance index of the system is given
by $P_{0}(R,d),$ and we wish to ensure that the performance is dependent
only on the evolution of observation and the initial state $Z_{\,0}.$

\subsection{The accounting context}

The dividend irrelevance question (in particular whether or not dividends
should be irrelevant to stockholders) has been a live issue since the 1961
paper \cite{MilM2} of Modigliani and Miller. See, for example, \cite{DybZ}. The
current quest for dividend irrelevance comes from the general possibility of
restating equity in terms of an identically discounted alternative series
based on accounting numbers, as first pointed out in 1936 by Preinreich
\cite{Pre}; see the discussion in the survey paper by \cite{OhlG}. If ($\Omega $)
models the evolution of the firm and $z_{t}$ models its observable
accounting numbers, interest focuses on whether valuations are possible at
time $t=0,$ based on the accounting numbers alone, that is to say in the
absence of access to the currently unobservable information $\omega _{\text{%
div}}$. See Proposition 2 below.

The current paper inter alia identifies circumstances under which dividend
irrelevance does indeed occur at an eigenvalue of $A$ and so shows that the
earlier derived equivalence is non-vacuous.

\subsection{Organization of material}

The paper is organized as follows. In Section 2 we give our main theorems
(Theorems 1-5) and the auxiliary propositions on which they are based.
Shorter proofs are included here, but longer proofs are delayed till later.
The results of Sections 2 and 3 are then used in Section 4 to perform a
detailed study of eigenvalue location of the augmented matrix. This is done
by examining the two-pole case first, and then estimating the distortion
effects when other poles are present. Some bifurcation analysis is conducted
in Section 5 in circumstances corresponding to Theorem 1. Sections 6 and
onwards contain longer proofs, or such details as are not required for the
analysis of Section 4.

\section{Main Theorems and auxiliary propositions}

In \cite{GaoOO} it is shown that dividend irrelevance at $R$ occurs iff $R$ takes
the value of the dominant eigenvalue, here defined to be the largest in
modulus (in the spirit of the Perron-Frobenius context -- see \cite[Ch. 8]{HorJ},
\cite[Ch.1,2]{Sen}), of the reduced matrix $A,$ which will forthwith be diagonal.
Asymptotic considerations suggest this result, since for generic initial
conditions, and for large $t$ the dominant eigenvalue growth of $%
A^{\lrcorner }$ dwarfs into insignificance the other state components, both
those entering the accounting state vector and those entering the dividend
(provided of course that the dividend-policy vector gives the dominant
growth component a non-zero coefficient). Asymptotic considerations thus
turn the multi-dimensional system apparently into an essentially
one-dimensional one, and it is to this that Ohlson's Principle (initially
proven in dimension one only) might apply -- see Theorem 2 below. That is to
say, assuming dividend and dominant state variable are inter-linked,
dividend irrelevance occurs if and only if $R$ takes a unique value, that
value being the dominant eigenvalue of the dividend-retention matrix $A$
(that of the dominant state). (Of course, in the long run, observation of
the dividend sequence permits inference of the dividend-policy vector.)

In this paper we offer an analysis of the quoted result based on algebraic
considerations, some complex analysis (including an inessential reference to
Marden's `Mean-Value Theorem for polynomials'), and graphical analysis.
These complement a standard textbook analysis based on Gerschgorin's circle
theorem -- for which see e.g. \cite[Th. 13.14]{Nob}  or \cite[Th. 7.8d]{Hen}.

Unsurprisingly, \textit{the eigenvalues of }$A^{\lrcorner }$\textit{\ may be
located arbitrarily,} but only if no restrictions are placed on the dividend
policy $(w,\beta ).$ Evidently, Dividend Irrelevance must implicitly assume
the convergence assumption as a bound on the eigenvalues of $A^{\lrcorner }$. It transpires (see Proposition 3) that the dividend-policy vector is
restricted by this assumption to the interior of an appropriate polytope in $%
\mathbb{R}^{n+1}$.

Conditions may be placed on the vector $b$ such that, when $\omega =(w,\beta
)$ lies in an open region of parameter space, it is the case that the
dominant eigenvalue of the augmented matrix is real and lies \textit{between
the first largest and the second largest} eigenvalue of the reduced matrix.
This is the substance of our first main result stated here and proved in
Sections 4 and 6.

\bigskip

\noindent \textbf{Theorem 1 (An Eigenvalue Dominance Theorem).} \textit{%
Suppose that }$A$ \textit{has real positive distinct eigenvalues. In the
canonical setting (\ref{canon}) we have as follows.}\medskip

\noindent (i) \textit{If }sign$[\delta _{1}]=-1$\textit{\ and }sign$[\delta
_{j}]=+1$\textit{\ for }$j=2,...,n,$\textit{\ then the open set}%
\begin{equation*}
\{\omega :A^{\lrcorner }\textit{ has real distinct roots and }\lambda
_{2}<\kappa _{1}(\omega )<\lambda _{1}\},
\end{equation*}%
\textit{has non-empty intersection with the set }%
\begin{equation*}
\{\omega :\omega _{1}>0,...,\omega _{n+1}>0\}.
\end{equation*}%
\textit{Moreover, }$\kappa _{2},$ \textit{the second largest eigenvalue of }$%
A^{\lrcorner }$ \textit{for small }$\omega ,$\textit{\ is increasing in }$%
\omega $ \textit{for small }$\omega .$ \textit{Under these circumstances
dividend irrelevance holds uniquely at }$R=\lambda _{1}.$\medskip

\noindent (ii) \textit{More generally, the open set}%
\begin{equation*}
\{\omega :A^{\lrcorner }\text{ has real distinct roots and }\lambda
_{2}<\kappa _{1}(\omega )<\lambda _{1}\},
\end{equation*}%
\textit{has non-empty intersection with the set }%
\begin{equation*}
\{\omega :\delta _{1}\omega _{1}<0,\delta _{2}\omega _{2}>0,...,\delta
_{n}\omega _{n}>0\},
\end{equation*}%
\textit{and again under these circumstances dividend irrelevance holds
uniquely at }$R=\lambda _{1}.$\medskip

\noindent (iii)\textit{\ If }$\delta _{1}\omega _{1}<0$\textit{\ and }$%
\delta _{k}\omega _{k}>0$\textit{\ for all }$k=2,..,n,$\textit{\ and }$%
\Omega $\textit{\ has all its eigenvalues in the disc }$|\zeta |<\lambda
_{1} $\textit{\ of the complex }$\zeta $\textit{-plane, then }$A^{\lrcorner }$\textit{%
\ has an eigenvalue in the annulus }%
\begin{equation*}
\mathcal{A}:=\{\zeta \in \mathbb{C}:\lambda _{2}<|\zeta |<\lambda _{1}\}\mathit{.}
\end{equation*}

\noindent (iv)\textit{\ If }$\delta _{1}\omega _{1}<0$\textit{\ and }$\delta
_{k}\omega _{k}<0$\textit{\ for all }$k=2,..,n,$\textit{\ then the system }$%
\Omega $\textit{\ has a real eigenvalue in the real interval }$(\lambda
_{2},\lambda _{1})$\textit{.}

\bigskip

For a proof see Section 6.

\bigskip

\noindent \textit{Remark.} We see therefore that for an appropriate vector $%
b $ there is a region of parameter space for which the eigenvalues of the
\textit{augmented} matrix $A^{\lrcorner }$ remain strictly bounded in
modulus by $\lambda _{1}$, the dominant eigenvalue of $A.$ Note the
re-emergence of the side conditions $\delta _{1}\omega _{1}<0$ analogous to
the condition $\omega _{12}\omega _{21}<0$ in Ohlson's Theorem for $n=1$
(see \cite{OhlG}).

\bigskip

We are able to provide some information about the extent of the subspace
(see formula (\ref{omegabound}) of \S 3) where we obtain (when $\delta
_{1}<0 $) the upper bound on positive $\omega _{1}$ of%
\begin{equation*}
\frac{1}{4}(\lambda _{1}-\beta )^{2}+\{\omega _{2}\delta _{2}+...\},
\end{equation*}%
for the case $\delta _{2}\omega _{2}>0.$ Moreover, Proposition 4 and
calculations of \S 4 appear to imply that, even if $\omega _{1}$ rises above
this bound, the two particular roots of the characteristic polynomial of $%
A^{\lrcorner }$ which are forced into coincidence remain outside the disc $%
|\zeta |\leq \lambda _{2}$ in the complex $\zeta $-plane (as they move
asymptotically to a vertical towards $\Re(\zeta )=\lambda _{2}$),
provided
\begin{equation*}
\beta >2\lambda _{2}-\lambda _{1}.
\end{equation*}%
By contrast, we find for $\omega _{1}\delta _{1}<0$ and $\omega _{2}\delta
_{2}<0$ the top two roots of the augmented matrix $A^{\lrcorner }$ both
approach $\lambda _{2}$ from opposite sides; this again is in keeping with
the expectation that dividend irrelevance occurs only at the dominant root $%
\lambda _{1}.$

\bigskip

Our results link to work concerned with the real spectral radius of a
matrix, see Hinrichsen and Kelb \cite{HinK}, which investigates by how much a
matrix may be perturbed without moving its spectrum out of a given open set
in the complex plane. In the cited work the open set of concern is usually
either the unit disc or the open left half-plane, both in connection with
stability issues. Our interest, however, focuses additionally on the open
set described by the \textit{annulus} $\mathcal{A}$ defined by the first and
second largest eigenvalues of $A$ (cf. Th. 1). We note that there is a
well-established Sturmian algorithm for counting the number of zeros of a
polynomial in the unit disc in the complex plane (see Marden \cite[\S 42, p.
148]{Mar}), and so in principle the issue of Dividend Irrelevance is resolvable
for a given policy vector $\omega _{\text{div}}$ by reference to the number
of zeros in the unit circle of the two polynomials%
\begin{equation*}
\chi _{A^{\lrcorner }}(\kappa /\lambda _{1}^{A}),\qquad \chi _{A^{\lrcorner
}}(\kappa /\lambda _{2}^{A}).
\end{equation*}%
Specifically, the first should have $n+1$ zeros and the second no more than $%
n.$ The Schur-Cohn criterion \cite[Th. 43.1]{Mar}, \cite[\S 6.8]{Hen} might perhaps
also be invoked to count the number of roots in the unit disc.

\subsection{Preliminaries}

Our analysis is based on two results embodied in Proposition 1 and in the
equivalences given in Proposition 2. The arbitrary placement of the zeros,
the substance of Proposition 3, is also a consequence of Proposition 2.

\bigskip

\noindent \textbf{Proposition 1 (Under the assumption of distinct
eigenvalues). }\textit{In the canonical setting (\ref{canon}) with }%
\begin{equation*}
\delta _{j}=\pm 1:
\end{equation*}%
\textit{for any }$j\leq n$\textit{\ and fixed }$R,$ \textit{the equity} $%
P_{0}(R,\omega )$ \textit{is locally or globally independent of }$\omega $
\textit{iff }$R=\lambda _{j}$\textit{\ provided}%
\begin{equation}
R>\max \{|\kappa _{k}(\omega )|:k=1,...,n+1\},  \label{convergence}
\end{equation}%
\textit{in which case}%
\begin{equation*}
d_{0}+P_{0}(R;d)=-\frac{RZ_{0}^{j}}{\delta _{j}}.
\end{equation*}

The proof is in Section 7.

\bigskip

\noindent \textit{Remark 1.} Apparently, if the eigenvalues of $A^{\lrcorner
}$ all lie in the disc with radius any other eigenvalue of $A,$ the
Proposition permits \textit{local} dividend irrelevance to occur at several
rates of return. We will show below that subject to (\ref{convergence}) such
an anomalous behaviour is definitely excluded when $\omega _{1}\neq 0$ and
also $\omega _{j}\neq 0$ for some $1<j\leq n.$

\noindent \textit{Remark 2.}\textbf{\ }In principle we might want to allow $%
\delta _{j}=-1,$ to respect a restriction in the directional sense of a
re-scaling of accounting variables (if appropriate); it transpires from the
next Proposition that the sign of $\delta _{j}$ can be absorbed by $\omega
_{j}$ and the choice of sign is only a matter of symbolic convenience, so
that we can interpret $\delta _{1}\omega _{1}<0$ as saying $\omega _{1}>0.$
That said, it is important to realize that rescaling an accounting variable,
say $z_{k}$ by $\alpha ,$ requires an inverse rescaling of the corresponding
dividend-policy component, that is of $\omega _{k}$ by $\alpha ^{-1}$ (in
order to preserve the definition of dividend untouched). The right-hand side
of the valuation equation perforce does not refer to the eigenvalues $\kappa
_{j}$, despite the fact that these control the growth rates of the canonical
accounting variables.

The following algebraic equivalences lie at the heart of all our arguments.
Below we denote by $\chi _{A^{\lrcorner }}(\kappa ,\omega _{1},...,\omega
_{n+1})$ the characteristic polynomial of $A^{\lrcorner }(\omega
_{1},...,\omega _{n+1}).$

\bigskip

\noindent \textbf{Proposition 2 (Inverse relations).} \textit{Put }$\lambda
_{n+1}=\beta $\textit{\ }$=\omega _{n+1}.$ \textit{The equations below are
all equivalent.}%
\begin{equation}
\chi _{A^{\lrcorner }}(\kappa ,\omega _{1},...,\omega _{n+1})=0.
\label{char}
\end{equation}%
\begin{align}\label{poly}
\prod\limits_{j=1}^{n+1}(\kappa -\lambda _{j}) =&\;\omega _{1}\delta
_{1}\prod\limits_{j=2}^{n}(\kappa -\lambda _{j})+\omega _{2}\delta
_{2}\prod\limits_{j\neq 2}^{n}(\kappa -\lambda _{j})\\
&\;+...+\omega _{n}\delta
_{n}\prod\limits_{j\neq n}^{n}(\kappa -\lambda _{j}).  \notag
\end{align}%
\textit{Polar form for} $j\leq n$ \textit{with \textbf{leading quadratic term%
}:}
\begin{equation}
\omega _{1}\delta _{1}+\omega _{2}\delta _{2}+...+\omega _{n}\delta
_{n}=(\kappa -\lambda _{j})(\kappa -\beta )+\sum_{k\neq j}^{n}\frac{\omega
_{k}\delta _{k}(\lambda _{j}-\lambda _{k})}{\kappa -\lambda _{k}}.\text{ }
\label{omega1}
\end{equation}%
\textit{Polar form for }$j=n+1$\textit{\ and with} $\kappa \neq \lambda _{k}$%
\textit{\ for }$k=1,...,n$\textit{:}
\begin{equation}
\beta =\kappa -\frac{\omega _{1}\delta _{1}}{\kappa -\lambda _{1}}-\frac{%
\omega _{2}\delta _{2}}{\kappa -\lambda _{2}}-...-\frac{\omega _{n}\delta
_{n}}{\kappa -\lambda _{n}}.  \label{beta}
\end{equation}%
\textit{In particular, with }$j=1,$ \textit{putting}%
\begin{equation}
f(\kappa ):=\{\omega _{2}\delta _{2}+..\}-(\kappa -\lambda _{1})(\kappa
-\beta )-\frac{\omega _{2}\delta _{2}(\lambda _{1}-\lambda _{2})}{\kappa
-\lambda _{2}}-...  \label{f-kappa}
\end{equation}%
\textit{we obtain the equivalent equation }%
\begin{equation*}
f(\kappa )=-\omega _{1}\delta _{1}.
\end{equation*}

Proof of equivalence follows in Section 8. Each of the above identities
enables a different analytic approach.

Our first conclusion regards the potentially arbitrary placement of the
zeros of (\ref{char}).

\bigskip

\noindent \textbf{Proposition 3 (Zero placement).} \textit{In the canonical
setting of Proposition 1, for an appropriate choice of real vector}$\omega $%
\textit{\ the characteristic polynomial}%
\begin{equation*}
\chi _{A^{\lrcorner }}(\kappa ,\omega )=|\kappa I-A^{\lrcorner }(\omega )|
\end{equation*}%
\textit{may take the form}%
\begin{equation}
\kappa ^{n+1}-p_{0}\kappa ^{n}+p_{2}\kappa ^{n-1}+...+(-1)^{n+1}p_{n},
\label{poly1}
\end{equation}%
\textit{for arbitrary choice of real coefficients }$p_{0},...,p_{n}.$ \textit{%
The transformation }$$(p_{0},...,p_{n})\rightarrow (\omega _{1},...,\omega
_{n+1})$$\textit{\ is affine invertible. The roots }$\kappa _{1},...,\kappa
_{n+1}$ \textit{of the characteristic polynomial may therefore be located at
will, subject only to the inclusion, for each selected complex root, of its
conjugate. }

\bigskip

This result is proved in Section 9.

Proposition 3 above indicates that in principle the region of parameter
space in which the boundedness assumption holds may be obtained as the
transform under the above mentioned transformation of the set of vectors $%
(p_{0},...,p_{n})$ satisfying a criterion derived from Cauchy's theorem on
the Inclusion Radius \cite[Th. 6.41]{Hen}, \cite[Th. 27.1]{Mar}, namely
\begin{equation*}
|p_{n}|+|p_{n-1}|\lambda _{1}+...+|p_{0}|\lambda _{1}^{n}<\lambda _{1}^{n+1}.
\end{equation*}%
(Recall that the inclusion radius of the polynomial (\ref{poly1}) is the
positive root of the polynomial $|p_{n}|+|p_{n-1}|\kappa +...+|p_{0}|\kappa
^{n}-\kappa ^{n+1}$.) Since the set of vectors $(p_{0},...,p_{n})$ so
described is the interior of a polytope, the corresponding region in
parameter space is therefore likewise seen to be the interior of a polytope.
Let us term this the \textbf{Cauchy polytope}.

Evidently $(0,...,0,\beta )$ is on the boundary of the Cauchy polytope,
since then%
\begin{equation*}
\chi _{A^{\lrcorner }}(\kappa ,\omega )=(\kappa -\lambda _{1})..(\kappa
-\lambda _{n})(\kappa -\beta ).
\end{equation*}

An immediate corollary is the following result, first announced for the case
$n=1$ by Ohlson at the 2003 International Conference on Advances in
Accounting-based Valuation -- see \cite[Lemma 4.1; generalization of Lemma
4.1: Appendix 2]{OhlG}.

\bigskip

\noindent \textbf{Theorem 2 (Multivariate Ohlson Principle).} \textit{The
system }$\Omega $\textit{\ has dividend irrelevance at }$R$\textit{\ iff }$%
R=\lambda _{1}.$

\medskip

\begin{proof} By varying $\omega $ we can place one eigenvalue
in the interval $(\lambda _{2},\lambda _{1}),$ so by Proposition 1, there
cannot be dividend irrelevance at $\lambda _{2}$ and below. Note that this
means that\textit{\ }$\delta _{1}\omega _{1}\neq 0$ for the chosen $\omega .$
\qed
\end{proof}

The situation with general placement of eigenvalues alters if $\beta $ is a
positive real, lies below the eigenvalues of $A,$ and the dividend-policy
vector $\omega $ of the canonical setting is non-negative in all its
components. The formula (\ref{beta}) confines the non-real eigenvalues $%
\kappa _{j}$ to an infinite strip, while the formula (\ref{poly}) allows us
to confine all the eigenvalues still further when $\omega $ is itself
bounded.

We refer to formula (\ref{omega1}) as the \textbf{associated polar form}.
This form offers a graphical approach to the analysis of the real root
location, and some insight into complex root location; in particular, the
\textbf{leading quadratic term} is responsible for unbounded root behaviour,
as follows.

\bigskip

\noindent \textbf{Proposition 4 (Unbounded roots).} \textit{Fix }$\omega
_{k} $\textit{\ for }$k\neq j$\textit{\ with }$$A_{j}=\sum_{h\neq j}\omega
_{h}\delta _{h}(\lambda _{j}-\lambda _{h})\neq 0.$$

\noindent (i) \textit{Subject to }$\lambda _{k}<|\kappa |$ \textit{we have
the asymptotic expansion}%
\begin{align*}
\omega _{1}\delta _{1}+&\;\omega _{2}\delta _{2}+...+\omega _{n}\delta
_{n}\\
=&\;(\kappa -\lambda _{j})(\kappa -\beta )+\sum_{s=1}^{\infty }\frac{1}{%
\kappa ^{s}}\left( \sum_{h\neq j}^{n}\omega _{h}\delta _{h}\lambda
_{h}^{s-1}(\lambda _{j}-\lambda _{h})\right) .
\end{align*}%
\textit{\noindent }(ii)\textit{\ For }$j\leq n$ \textit{the unbounded roots
as }$\delta _{j}\omega _{j}\rightarrow -\infty $\textit{\ behave
asymptotically as follows:}%
\begin{equation*}
\kappa =\left\{
\begin{array}{cc}
\pm \sqrt{|\omega _{j}|}+O(\omega _{j}^{-1/2}), & \text{for }\delta _{j}=1,
\\
\left[ \frac{1}{2}(\lambda _{j}+\beta )+\frac{A_{j}}{2|\omega _{j}|}%
++O(\omega _{j}^{-2})\right] \pm i\left[ \sqrt{|\omega _{j}|}+O(\omega
_{j}^{-1/2})\right] , & \text{for }\delta _{j}=-1.%
\end{array}%
\right.
\end{equation*}

For the proof, see Section 10.

See Figure $\ref{gig1}$ above and Figure $\ref{gig4}$ in Section 4 for illustrative examples.

\bigskip

\noindent \textit{Remark 1.} In the case $j=1$ with $\delta _{1}=-1,$ we are
of course assuming that $\omega _{1}\rightarrow \infty .$ If moreover $%
\delta _{k}\omega _{k}>0$ for all $k=2,...,n,$ and $\lambda _{1}>\lambda
_{2}>...>\lambda _{n},$ we have $A_{1}=\sum_{h\neq 1}\omega _{h}\delta
_{h}(\lambda _{1}-\lambda _{h})>0.$ Here the conjugate roots have real part
approaching $\frac{1}{2}(\lambda _{j}+\beta )$ from the right. However, with
other sign assumptions on $\delta _{h}\omega _{h},$ the sign of $A_{1}$ need
not be positive, in particular if $\delta _{h}\omega _{h}<0$ for all $h.$

\bigskip

\noindent \textit{Remark 2.}\textbf{\ }By (\ref{poly}) we may rewrite the
characteristic polynomial in the form
\begin{equation*}
\frac{1}{\omega _{j}}\prod\limits_{k=1}^{n+1}(\kappa -\lambda _{k})-\frac{1%
}{\omega _{j}}\{\omega _{1}\delta _{1}\prod\limits_{k=2}^{n}(\kappa
-\lambda _{k})+\omega _{2}\delta _{2}\prod\limits_{k\neq 2}^{n}(\kappa
-\lambda _{k})+...+\omega _{n}\delta _{n}\prod\limits_{k\neq n}^{n}(\kappa
-\lambda _{k})\}.
\end{equation*}%
For fixed\textit{\ }$\omega _{k},$\textit{\ }with\textit{\ }$j\neq k,$%
\textit{\ }pass to the limit as $|\omega _{j}|\rightarrow \infty $ to obtain
the following equation of degree $n-1:$%
\begin{equation*}
\prod\limits_{k\neq j}^{n}(\kappa -\lambda _{k})=0.
\end{equation*}%
Thus given the assumptions of the Proposition, \textbf{only two} complex
roots can be unbounded.

\bigskip

\noindent \textit{Remark 3.} Note that, by contrast, the unbounded roots for
$\omega _{1},...,\omega _{n}$ fixed and $\beta $ varying have the asymptotic
behaviour $\kappa =\beta +O(\beta ^{-1}).$ Note also that, if $\sum_{k\neq
j}\omega _{k}\delta _{k}(\lambda _{j}-\lambda _{k})=0,$ then the error term $%
O$-behaviour alters.

\bigskip

\noindent \textbf{Theorem 3}. \textit{With fixed }$\omega _{k}$\textit{\ for
}$k\neq 1$ \textit{such that }$\sum_{k\neq 1}\omega _{k}\delta _{k}(\lambda
_{j}-\lambda _{k})\neq 0,$ \textit{and with }$\delta _{1}=-1$\textit{\ if }%
\begin{equation*}
\beta \geq 2\lambda _{2}-\lambda _{1},
\end{equation*}%
\textit{\ the unbounded root locus does not enter the disc }$|\zeta |\leq
\lambda _{2}$\textit{\ as }$\omega _{1}\rightarrow \infty .$ \textit{So the
system }$\Omega $\textit{\ has local dividend irrelevance at }$(\omega
,\beta )$\textit{\ uniquely at} $R=\lambda _{1}.$

\bigskip

\begin{proof} Under these circumstances the unbounded
roots are outside the disc $|\zeta |\leq \lambda _{2}$, since they are
confined to $\Re(\zeta )>\lambda _{2}$ by virtue of
\begin{equation*}
\lambda _{2}\leq \frac{1}{2}(\lambda _{1}+\beta ).
\end{equation*}%
By Proposition 1 the only value remaining for $R$ is thus $\lambda _{1}.$\qed
\end{proof}

\bigskip

\noindent \textbf{Notation. }Below and throughout, $K(\varepsilon )$\textit{%
\ }denotes\textit{\ }the real interval\textit{\ }$[\beta -\eta (\varepsilon
),\lambda _{1}+\eta (\varepsilon )],$\textit{\ }where%
\begin{equation*}
\beta -\eta =\frac{(\beta +\lambda _{1})-\sqrt{(\lambda _{1}-\beta
)^{2}+4\varepsilon }}{2},\qquad \lambda _{1}+\eta =\frac{(\beta +\lambda
_{1})+\sqrt{(\lambda _{1}-\beta )^{2}+4\varepsilon }}{2};
\end{equation*}%
$S(K,\pi /(n+1))$\ comprises the two circles in the plane subtending angles
of $\pi /(n+1)$\ on $K(\varepsilon ).$

\bigskip

\noindent \textbf{Proposition 5 (Strip-and-two-circles theorem). } \textit{%
Suppose that }$\beta \leq \lambda _{n}<...<\lambda _{1},$ \textit{that }$%
\omega _{1}\neq 0,$ \textit{and that }%
\begin{equation*}
\delta _{1}\omega _{1},...,\delta _{n}\omega _{n}\geq 0.
\end{equation*}

\noindent (i) \textit{All the non-real roots of the characteristic equation (%
\ref{char}) lie in the infinite strip of the complex }$\zeta $\textit{-plane
given by}%
\begin{equation*}
\beta \leq \Re(\zeta )\leq \lambda _{1}.
\end{equation*}

\noindent (ii) \textit{For }$\varepsilon >0$\textit{\ arbitrary, if }%
\begin{equation*}
\omega _{1}+...+\omega _{n}\leq \varepsilon ,\qquad \delta _{1}\omega
_{1},...,\delta _{n}\omega _{n}\geq 0,
\end{equation*}%
\textit{then all the roots of (\ref{char}) lie in the star-shaped region }$%
S(K,\pi /(n+1))$\textit{.}

The proof is delayed to Section 11.

\bigskip

\noindent \textit{Remark 1.} Taken together parts (i) and (ii) may operate
simultaneously. These results should, however, be taken together with
Gerschgorin's Circle Theorem, which implies immediately that the eigenvalues
lie in the union of the discs in the complex $\zeta $-plane given by $|\zeta
-\lambda _{j}|\leq |\omega _{j}|$ and by $|\zeta -\beta |\leq |\omega
_{1}|+..|\omega _{n}|.$ Thus the eigenvalues are bounded, not only to the
above mentioned vertical strip but also to a horizontal strip of width $%
2\max \{|\omega _{j}|:$ $j\leq n\}$ around the real axis.

\bigskip

\noindent \textit{Remark 2.} It is obvious that, for $\omega _{2}=..=\omega
_{n}=0$ and with $|\omega _{1}|\leq \varepsilon ,$ the real roots of (\ref%
{char}) lie in $K(\varepsilon )$ by continuity. Gerschgorin's Circle Theorem
limits the real roots to the slightly larger interval $[\beta -\varepsilon
,\lambda _{1}+\varepsilon ].$ Thus the two-circle result is merely a
sharpening of the bounds.

\bigskip

\noindent \textit{Remark 3.} If $\lambda _{n}<\beta ,$ less elegant
improvements can be made so that $K$ extends only as far as $\lambda _{1}$
on the left.

We can state, ahead of the proof of Proposition 5, our theorem on eigenvalue
location.

\noindent \textbf{Theorem 4 (Eigenvalue bounds)}. \textit{Suppose that }$%
\beta \leq \lambda _{n}<...<\lambda _{1}$ \textit{and that }%
\begin{equation*}
|\omega _{1}|,...,|\omega _{n}|\leq \varepsilon .
\end{equation*}%
\textit{Non-real eigenvalues lie in the rectangle bounded by }$\gamma =\beta
,\zeta =\lambda _{1},\zeta =\pm \varepsilon .$ \textit{Real eigenvalues lie
in the interval }$K(\varepsilon ).$

\bigskip

\begin{figure}[!t] %F1
\centering\includegraphics[width=0.7\textwidth] {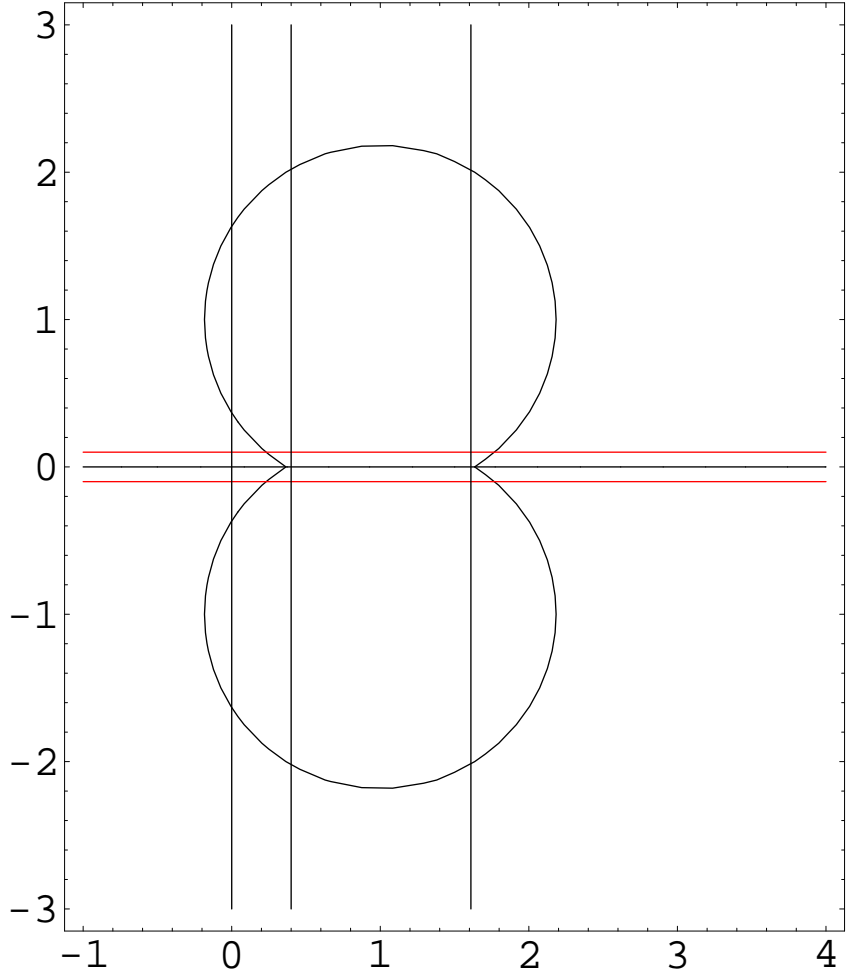}
\caption{
Vertical strip and two-circle bounds. Horizontal bound implied by Gerschgorin's Theorem.
} \label{gig2}
\end{figure}

The theorem follows from Proposition 5 - see Figure $\ref{gig2}$. The two-circle result
gives useful bounds only for the real roots.

\noindent \textit{Remark.} The above analysis does not yet exclude the
possibility of all eigenvalues being located to the left of $\lambda _{2}.$
We next offer a graphical analysis of the real-root locations in the
following subsection, which shows that at least one root has to be to the
right of $\lambda _{1}$ when $\omega _{1}\neq 0$ and $\omega _{j}\neq 0$ for
some $j>1.$

\section{Eigenvalue location: general analysis}

Our ultimate purpose (achieved in the next section) is to show that under
suitable restrictions one can guarantee the existence of a real eigenvalue
in the range $(\lambda _{2},\lambda _{1}).$ Specifically, we show that if $%
\delta _{1}=-1$ and $\delta _{2}=...=+1$ there is a real eigenvalue in the
range $(\lambda _{2},\lambda _{1})$ for all small enough positive $\omega
_{1}$. The aim of this section is to identify (i) the real-root locus of the
characteristic polynomial of $H(\omega )$ by examining the general features
of the graph of the associated polar form of the characteristic equation
given by (\ref{omega1}), and (ii) the complex-root locus when the associated
polar form has just two poles. The latter is a preliminary to our
identification of (iii) an elementary estimate of the distortion effect of
other poles.

In this section the real eigenvalues of $H(\omega )$ (see equation(\ref%
{canon})) are studied as functions of $\omega _{1}$ with the other
components of $\omega $ fixed. Interest naturally focuses on $\omega _{1}$
as the link coefficient with the dominant state vector. The two-pole case
arises when $n=3$ and is considered as a benchmark, with a view to
understanding how the multi-pole situation deforms the benchmark case.

Treating $\omega _{1}$ as a free variable, with the remaining
dividend-policy coefficient fixed, we use (\ref{omega1}) to study the map $%
\kappa \rightarrow \omega _{1}$ and its local inverses. We have, with $f$ as
in (\ref{f-kappa}),
\begin{equation}
-\omega _{1}\delta _{1}=f(\kappa )=\{\omega _{2}\delta _{2}+..\}-(\kappa
-\lambda _{1})(\kappa -\beta )-\frac{\omega _{2}\delta _{2}(\lambda
_{1}-\lambda _{2})}{\kappa -\lambda _{2}}-...,  \label{eqn}
\end{equation}%
so that the graph of $f$, or of $\omega _{1},$ against $\kappa $ has $(n-1)$
vertical asymptotes from right to left at $\kappa =\lambda _{2},\lambda
_{3},...\lambda _{n}$ all of which are manifestly simple poles. The
asymptotes break up the concave leading quadratic term (if $\delta _{1}<0)$
into $n$ connected components corresponding to the intervals $(-\infty
,\lambda _{n}),(\lambda _{n},\lambda _{n-1}),...,(\lambda _{n},+\infty )$.
The equation%
\begin{equation*}
\frac{\partial \omega _{1}}{\partial \kappa }=0
\end{equation*}%
is equivalent to an $n$-degree polynomial equation so its $n$ roots
contribute to at most $n$ stationary points in the graph.

In the interval $(\lambda _{j+1},\lambda _{j})$ the component has an even,
respectively an odd, number of stationary points depending on whether the
sign of $\omega _{j+1}\delta _{j+1}\omega _{j}\delta _{j}$ is $+1$ or $-1.$
In view of the behaviour of the leading quadratic term, not all the
components can be monotone (possess a zero number of stationary points!).
Thus at least one component is non-monotonic.

The components may be interpreted as graphs/loci of the eigenvalues $\kappa
_{j}(\omega _{1}).$ More precisely, the differentiable local inverses of the
mapping $\kappa \rightarrow \omega _{1}$ are the graphs of $\kappa
_{j}(\omega _{1}).$ That is to say, each non-monotonic component must be
first partitioned into monotone parts on either side of its stationary
points. The labelling of these inverses from right to left respects the
cyclic order on the set $\{1,...,n\}$ together with one or other of the
identifications%
\begin{equation*}
\lim_{\kappa \nearrow \lambda _{j}}\omega _{1}(\kappa )=\kappa _{j},\qquad
\lim_{\kappa \searrow \lambda _{j}}\omega _{1}(\kappa )=\kappa _{j}.
\end{equation*}%
The latter may require the point at infinity on the asymptote $\kappa
=\lambda _{j}$ to be considered as the intersection of consecutive loci.

Note that from (\ref{eqn}) $f(\lambda _{1})=0$ and so $\kappa
_{1}(0)=\lambda _{1}$. (This is consistent with the matrix $H(\omega
)-\kappa I$~ having a first column with zeros in all but the last row.)

We will see in \S 4.1 from these asymptotic features of the graph that,
since $\kappa _{1}(0)=\lambda _{1},$ for all small enough positive $\omega
_{1}$ the eigenvalue $\kappa _{1}(\omega _{1})$ is in the range $(\lambda
_{2},\lambda _{1}),$ as we now demonstrate.

As a preview of the full argument of \S 4.2, with our assumption that $%
\delta _{1}=-1$ and $\delta _{2}=\delta _{3}=...=+1,$ we note that we can
arrange for $\kappa _{1}(\omega _{1})$ to be large and positive in the
vicinity to the right of $\lambda _{1}$ by taking $\omega _{2}<0.$ With $%
\omega _{2}<0$ the domain of $\kappa _{1}$ is infinite, so that%
\begin{equation*}
\lim_{\omega _{1}\rightarrow \infty }\kappa _{1}(\omega _{1})=\lambda _{2}.
\end{equation*}

Thus the largest real eigenvalue of $\kappa _{1}$ remains above $\lambda
_{2}.$ See Figure $\ref{gig3a}$ (in \S 4.1 below). Of course for small enough $\omega
_{1}$ the remaining roots $\kappa _{j}(\omega _{1}),$ even if complex,
remain in an open vertical complex strip including the closed real interval $%
[\beta ,\lambda _{2}].$

We can similarly arrange for $\kappa _{1}(\omega _{1})$ to be large and
negative in the vicinity to the right of $\lambda _{1}$ by taking $\omega
_{2}>0.$ In view of the behaviour of the graph for large $\kappa >\lambda
_{1},$ this implies the existence of two roots in $(\lambda _{1},\lambda
_{2})$ under these circumstances. With $\omega _{2}>0$ the domain of $\kappa
_{1}$ is bounded, say by $\omega _{1}\leq \omega _{1}^{\ast }=\omega
_{1}^{\ast }(\omega _{2},...,\omega _{n+1}),$ and one has
\begin{equation*}
\lim_{\omega _{1}\searrow \omega _{1}^{\ast }}\kappa _{1}(\omega
_{1})=\lim_{\omega _{1}\nearrow \omega _{1}^{\ast }}\kappa _{2}(\omega _{1}).
\end{equation*}%
See Figure $\ref{gig3b}$ (below in \S 4.1). As the eigenvalue $\kappa _{1}$ remains
above $\lambda _{2},$ dividend irrelevance can occur only at $\lambda _{1}.$
An upper bound for $\omega _{1}^{\ast }$ is provided by the maximum value of
the leading quadratic term%
\begin{equation*}
-(\kappa -\lambda _{1})(\kappa -\beta )+\{\omega _{2}\delta _{2}+..\}
\end{equation*}%
obtained by evaluation at $\kappa =\frac{1}{2}(\lambda _{1}+\beta ),$ namely%
\begin{equation}
\frac{1}{4}(\lambda _{1}-\beta )^{2}+\{\omega _{2}\delta _{2}+...\}.
\label{omegabound}
\end{equation}%
This gives $\beta ,$ the coefficient at the previous date's dividend, a
significant bounding role.

In order to understand the qualitative behaviour of the complex eigenvalues
of the characteristic equation, we study the associated polar form on an
open interval between the two adjacent poles at $\lambda _{j+1}$and $\lambda
_{j}$ for $j\geq 2.$ As a first step, we study the contribution to $f(\kappa
)$ arising in (\ref{eqn}) only from the two terms corresponding to the two
adjacent poles $\lambda _{j+1}$ and $\lambda _{j}$. In the subsequent
section we identify how the presence of the other terms in (\ref{eqn})
perturbs this simple analysis.

\subsection{Root locus for two-pole case}

For present purposes, by scaling and a shift of origin, as $\lambda
_{j+1}>\lambda _{j}$ we may take $\lambda _{j+1}=0$ and $\lambda _{j}=\tau
>0.$ It is convenient to study the corresponding terms of $f(\kappa )$ by
introducing $T>0$ and looking at the two functions%
\begin{eqnarray*}
f_{1}(x) &=&\frac{1}{x}+\frac{T}{x-\tau }=\frac{(1+T)x-\tau }{x(x-\tau )}, \\
f_{2}(x) &=&\frac{1}{x}-\frac{T^{2}}{x-\tau }.
\end{eqnarray*}%
The first of these has its zero at $x=\tau /(1+T)$ in $(0,\tau )$ and maps $%
(0,\tau )$ bijectively onto the reals. The second function is more awkward;
provided $T\neq 1,$ it has a zero outside $(0,\tau )$ at $\tau /(1-T^{2}).$
More information is provided in the Circle Lemma below.

\bigskip

\noindent \textbf{Circle Lemma.} \textit{For }$0<T\neq 1,$\textit{\ the
function }$f_{2}(x)$\textit{\ has a positive local minimum value }$K_{+}$%
\textit{\ at }$x^{+},$\textit{\ and positive local maximum value }$K_{-}$%
\textit{\ at }$x_{-},$\textit{\ where}%
\begin{equation*}
x^{\pm }=\frac{\tau }{1\pm T},\qquad K_{\pm }=\frac{(1\pm T)^{2}}{\tau }.
\end{equation*}%
\textit{The range of }$f_{2}(x)$\textit{\ on }$(0,\tau )$\textit{\ omits an
interval of positive values }$(K_{-},K_{+})$.

\textit{For }$K_{-}<K<K_{+}$\textit{\ the equation}%
\begin{equation}
f_{2}(z)=K,\text{ equivalently }z^{2}-(\tau +(1-T^{2})K^{-1})z+\tau K^{-1}=0,
\label{quad}
\end{equation}%
\textit{has conjugate complex roots in the complex }$\zeta $\textit{-plane
describing a circle centred at the real number }$(x^{+}+x^{-})/2$\textit{\
with radius }$|x^{+}-x^{-}|/2.$

\textit{The real part moves from }$x^{+}$\textit{\ towards }$x^{-}.$ \textit{%
Hence for }$0<T<1,$\textit{\ as }$K$ \textit{decreases from }$K_{+}$\textit{%
\ the real part increases and for }$T>1$ \textit{it decreases.}

\bigskip

\textit{Note.} See below for a more general analysis of the behaviour of the
real part near a local minimum. For $T=1$ the function is symmetric about $%
x=\tau /2$ and is asymptotic to zero at infinity; a limiting version of the
lemma is thus still valid, but the conjugate roots lie on the vertical line $%
\Re(\zeta )=\tau /2$ for $0<K<K_{+}=4/\tau .$

\bigskip

\begin{proof} The conjugate roots $z,\bar{z}=x\pm iy$ satisfy
\begin{equation*}
x^{2}+y^{2}=z\bar{z}=\tau K^{-1}.
\end{equation*}%
In view of the dependence of the real part on $K$ as given by $x=$ $(\tau
+(1-T^{2})K^{-1})/2,$ we see that the term $\tau K^{-1}$ can be absorbed by
a shift of origin on the $x$ axis. Hence the locus as $K$ varies is a
circle. Reference to the extreme locations $x^{\pm }$ identifies the shifted
centre, and the radius.\qed
\end{proof}

\subsection{Multi-pole case: a distortion estimate}

The presence of additional poles outside $(\lambda _{j+1},\lambda _{j})$
distorts the result obtained in the Circle Lemma above for the two-pole
contribution. Provided all the eigenvalues $\lambda _{j}$ are
well-separated, i.e. the ratio of adjacent intervals does not vary greatly
(see the calculation of $\tau /c$ below), the distortion is controlled by
the value of $\sum_{k\notin \{j,j+1\}}\omega _{k}\delta _{k}(\lambda
_{j}-\lambda _{k})$. In the next subsection we take note of a
third-derivative test (based on Taylor's Theorem) which identifies
bifurcation behaviour of coincident real roots. For an application see \S %
4.5.1.

We again work in the standardized co-ordinate system with the adjacent poles
at $x=0$ and $x=\tau .$ Put $\xi :=\tau /2;$ then the decomposition, valid
for $x\in (0,\tau )$ and $c>0,$
\begin{equation*}
\frac{1}{x+c}=\frac{A_{0}(x)}{x}1_{[\xi ,a]}(x)-\frac{A_{\tau }(x)}{x-\tau }%
1_{[0,\xi ]}(x)
\end{equation*}%
yields the following bounds:%
\begin{eqnarray*}
0 &<&\frac{\tau }{2c+\tau }<A_{0}(x)<\frac{c}{c+\tau }<1,\text{ for }\xi
<x<a, \\
0 &<&\frac{\tau }{2c+\tau }<A_{\tau }(x)<\frac{\tau }{c},\text{ for }0<x<\xi
.
\end{eqnarray*}%
These and an amendment of the parameter $T$ in $f_{2}(x)$ to a variable
coefficient $T=T(x)$ for $x$ in the closed interval $[0,\tau ]$ enable us to
account for the presence of terms other than $f_{2}(x)$ in $f(x),$ by
absorbing their contributions into $T(x).$ On $[0,\tau ],$ the `adjusted' $T$
inherits from $f$ boundedness, continuity and indeed two-fold
differentiability. Thus we have%
\begin{equation*}
f(x)=\frac{1}{x}-\frac{T}{x-\tau },\qquad f^{\prime }=-\frac{1}{x^{2}}+\frac{%
S}{(x-\tau )^{2}},
\end{equation*}%
where $S=S(x)=T(x)-(x-\tau )T^{\prime }(x)\sim T(\tau ),$ for $x$ close to $%
\tau $. Expansion round $x$ yields
\begin{equation*}
T(\tau )=T(x)+(\tau -x)T^{\prime }(\tau )+\frac{1}{2}(\tau -x)^{2}T^{\prime
\prime }(x)+o((\tau -x)^{3}),\qquad \text{as }x\rightarrow \tau ,
\end{equation*}%
so%
\begin{equation*}
S=T(\tau )-\frac{1}{2}(x-\tau )^{2}T^{\prime \prime }(x)+o((\tau -x)^{3}).
\end{equation*}

The main point of this is to observe how $S$ perturbs the complex roots (cf.
Figure $\ref{gig5b}$). The equation (\ref{quad}) of the Circle Lemma now gives us the
following.

\bigskip

\noindent \textbf{Proposition 6 (Distortion Estimate).} \textit{The equation}%
\begin{equation*}
K=\frac{1}{x}-\frac{T(x)}{x-\tau }
\end{equation*}%
\textit{is equivalent to }$f_{2}(x)=K+\bar{T}$\textit{\ for}%
\begin{equation*}
\bar{T}=T^{\prime }(x)+\frac{1}{2}(\tau -x)T^{\prime \prime }(x)+o((\tau
-x)^{2}).
\end{equation*}%
\textit{The equation }$f(x)=K\ $\textit{is equivalent to}%
\begin{equation*}
x=\frac{\tau }{2}+\frac{(1-S^{2})}{2K}\pm \frac{i}{2}\sqrt{\delta (K^{-1})},
\end{equation*}%
\textit{with }$\delta (k)=\tau ^{2}K_{+}K_{-}(k-k_{+})(k_{-}-k),$ \textit{%
where}
\begin{eqnarray*}
0 &<&k_{-}^{-1}=K_{-}=\frac{(1-S(x))^{2}}{\tau }<K<\frac{(1+S(x))^{2}}{\tau }%
=K_{+}=k_{+}^{-1}, \\
S &=&T(\tau )-\frac{1}{2}(x-\tau )^{2}T^{\prime \prime }(x)+o((\tau -x)^{3}.
\end{eqnarray*}

\bigskip

We note that $\frac{1}{2}T^{\prime \prime }(x^{+})=-f^{\prime
}(x^{+})+o((\tau -x^{+})).$

\subsection{Analysis of the real part via Taylor's theorem}

\noindent \textbf{Proposition 7 (Third-derivative test: real part follows }$%
f^{\prime \prime }/f^{\prime \prime \prime }$\textbf{).} \textit{Suppose
that }$f(\kappa )$ \textit{has a local minimum/maximum} \textit{at }$\kappa
=\kappa ^{\ast }.$ \textit{Let }$\kappa (\omega )$\textit{\ denote the local
solution for }$\kappa $ \textit{over the complex domain of the equation }$%
g(\kappa )=0$ \textit{for}%
\begin{equation*}
g(\kappa )=\left\{
\begin{array}{cc}
f(\kappa )-\omega ^{2}, & \text{ for }\kappa ^{\ast }\text{ a local minimum,}
\\
f(\kappa )+\omega ^{2}, & \text{for }\kappa ^{\ast }\text{ a local maximum,}%
\end{array}%
\right.
\end{equation*}%
\textit{with }$\omega >0$\textit{\ small and with }$\kappa (0)=\kappa ^{\ast
}$\textit{. If }$f^{\prime \prime \prime }(\kappa ^{\ast })\neq 0,$\textit{\
then the locus }$\kappa (\omega )$\textit{\ satisfies initially: }%
\begin{equation*}
\Re(\kappa (\omega ))\text{\textit{\ is increasing\ if}}%
\mathit{\ }f^{\prime \prime }(\kappa ^{\ast })/f^{\prime \prime \prime
}(\kappa ^{\ast })>0,\mathit{\ }\text{\textit{and} \textit{ is otherwise
decreasing}}.
\end{equation*}

\bigskip

For a proof see Section 12. Proposition 7 implies that the quadratic terms
of the associated polar form $f(\kappa )$ have no effect on the local
behaviour of the real part at a bifurcation.

\section{Eigenvalue location:\ some cases}

In this section we consider the case $n=2,$ and the two cases with $n$
general when $\delta _{k}\omega _{k}$ are of constant sign for $k=2,...,n,$
as referred to in Theorem 1.

\subsection{The case $n=2$}

This case is in fact typical, despite having the simplifying structure that
one root of the cubic characteristic polynomial $\chi _{H}$ is always real.
There may thus be two more real roots, or two conjugate complex roots.

In view of earlier comments, we need to consider only the case $\delta
_{1}\omega _{1}<0.$ Interpret this as saying $\delta _{1}=-1$ and $\omega
_{1}>0.$

\begin{figure}[!t] %F1
\centering\includegraphics[width=0.7\textwidth] {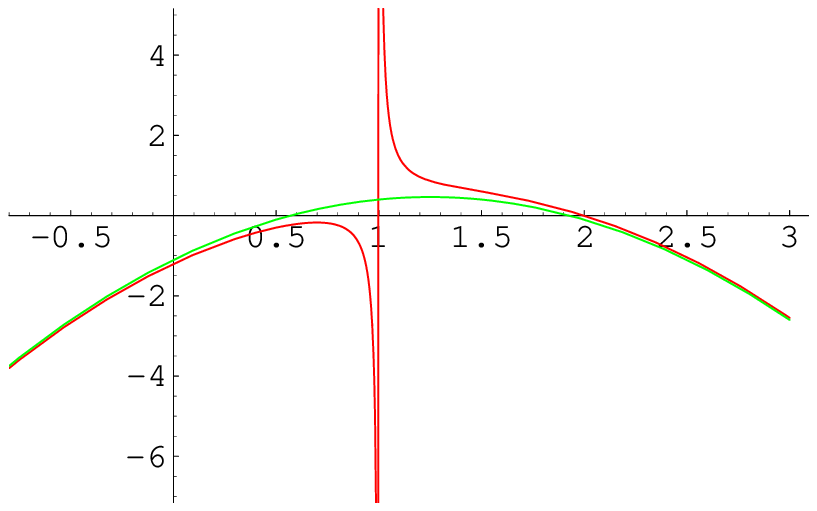}
\caption{
Graph of $\omega_1(\kappa)$ with $\delta_1=-1,\,\delta_2=+1,\,\omega_2<0$.
Leading quadratic in green; $\kappa$-axis horizontal.
} \label{gig3a}
\end{figure}

\begin{figure}[!t] %F1
\centering\includegraphics[width=0.7\textwidth] {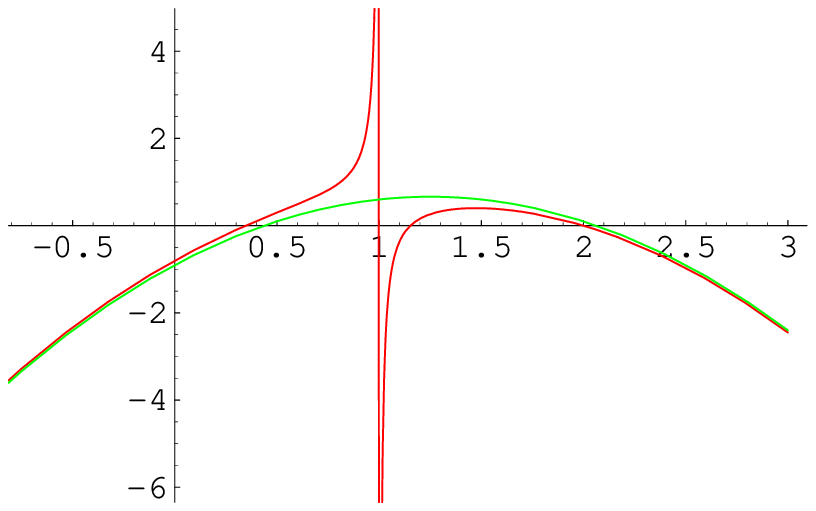}
\caption{
Graph of $\omega_1(\kappa)$ with $\delta_1=-1,\,\delta_2=+1,\,\omega_2>0$.
Leading quadratic in green; $\kappa$-axis horizontal.
} \label{gig3b}
\end{figure}

(a) For now assume $\delta _{2}\neq 0.$(For $\delta _{2}=0$ see (b) below.)
Taken together Figures $\ref{gig3a}$ and $\ref{gig3b}$ tell it all. They graph the implicit
relation between $\omega _{1}$ and the eigenvalues $\kappa $ as given by the
equation $\chi _{H}(\kappa ;\omega _{1})=0,$ treating $\kappa $ as the
independent variable and $\omega _{1}$ as dependent. To derive the root
locus for real roots $\kappa $ rotate the graphs, so that $\omega _{1}$
becomes the independent variable. Then each branch of the graph yields the $%
\kappa _{j}$ as the dependent eigenvalues in decreasing magnitude.%

If $\delta _{2}\omega _{2}<0,$ then for increasing $\omega _{1},$ as in
Figure $\ref{gig3a}$, the dominant root of $H(\omega )$ decreases down to $%
\lambda _{2}$ (in the limit).

If $\delta _{2}\omega _{2}>0,$ then for increasing $\omega _{1},$ as in
Figure $\ref{gig3b}$, the first/second root of $H(\omega )$,
respectively, decrease/increase into coincidence in the interval $(\lambda
_{2},\lambda _{1}).$ Thereafter, the root locus of the conjugate pair
behaves as illustrated below in Figure $\ref{gig4}$. That is, the real part decreases
towards $\frac{1}{2}(\lambda _{1}+\beta )$ and the imaginary parts tend to
infinity.%
\begin{figure}[!t] %F1
\centering\includegraphics[width=0.7\textwidth] {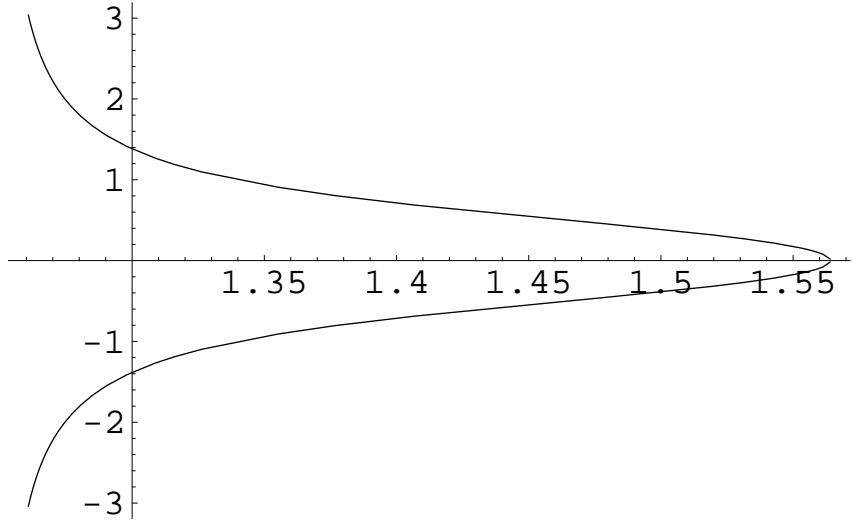}
\caption{
Locus of the conjugate roots in the complex $\zeta$-plane for $\delta_1<0$ as $\omega_1$ increases with $\omega_2$ fixed. Vertical asymptote is $\Re(\zeta)=\frac{1}{2}(\lambda_1+\beta).$
} \label{gig4}
\end{figure}

(b) We consider now the limiting case $\delta _{2}=0$. As the characteristic
polynomial is
\begin{eqnarray*}
\chi _{H}(\kappa ) &=&(\kappa -\lambda _{1})(\kappa -\lambda _{2})(\kappa
-\beta )-\omega _{1}\delta _{1}(\kappa -\lambda _{2}) \\
&=&(\kappa -\lambda _{2})[(\kappa -\lambda _{1})(\kappa -\beta )-\omega
_{1}\delta _{1}],
\end{eqnarray*}%
this has one root at $\kappa =\lambda _{2}$ which, by Proposition 1,
precludes DPI occurring at $R=\lambda _{2}.$ Note that the expression in
square brackets corresponds to an analogous problem with reduced order
(effectively the $n=1$ case), and this feature, of a reduction of order,
occurs for general $n$.

\bigskip

\noindent \textit{Remark.} Of course, for $\delta _{1}\omega _{1}<0,$ the
characteristic polynomial of $H(\omega )$ for $\delta _{2}=0$ has two
additional real roots in the interval $(\beta ,\lambda _{1})$ for small $%
|\omega _{1}|.$ For larger $|\omega _{1}|$ the conjugate roots lie in the
complex $\zeta $-plane on the vertical $\Re(\zeta )=\frac{1}{2}%
(\lambda _{1}+\beta ).$ This is because
\begin{equation*}
\kappa =\frac{(\beta +\lambda _{1})\pm \sqrt{(\lambda _{1}-\beta
)^{2}+4\omega _{1}\delta _{1}}}{2}=\frac{1}{2}(\lambda _{1}+\beta )\pm i\{%
\sqrt{|\omega _{1}|}+O(|\omega _{1}|^{-1/2})\}.
\end{equation*}

As before, if $\beta \geq 2\lambda _{2}-\lambda _{1}$\textit{, }then $%
\lambda _{2}\leq \frac{1}{2}(\lambda _{1}+\beta )$ and the conjugate roots
do not enter the open disc $|\zeta |<\lambda _{2}$; however, if $\beta
<2\lambda _{2}-\lambda _{1}$\textit{, }then $\frac{1}{2}(\lambda _{1}+\beta
)<\lambda _{2},$ which means that two roots will be in the open disc for a
range of values of $|\omega _{1}|$ with the third on the boundary.
Eventually the conjugate pair exits the annulus $\mathcal{A}$ (see Figure $\ref{gig1}$).

\subsection{General case $\protect\delta _{1}\protect\omega _{1}<0$ with $%
\protect\delta _{j}\protect\omega _{j}>0$ for $j=2,...,n.$}

We fix arbitrarily $\delta _{j}\omega _{j}>0$ for $j=2,...,n$. The
assumption $\omega _{1}\delta _{1}<0$ is without loss of generality
interpreted as $\delta _{1}=-1$ with variable $\omega _{1}>0.$ The analysis
now proceeds similarly. Here we have the identity connecting eigenvalue $%
\kappa $ and parameter $\omega _{1}$ in the shape of the associated polar
form (\ref{omega1}):
\begin{equation*}
\omega _{1}=f(\kappa ):=-(\kappa -\lambda _{1})(\kappa -\beta )+\{\omega
_{2}\delta _{2}+..\}-\frac{\omega _{2}\delta _{2}(\lambda _{1}-\lambda _{2})%
}{\kappa -\lambda _{2}}-....
\end{equation*}%
We note that the quadratic term reflects the behaviour of the two equations
which remain when all the variables other than the dividend and the dominant
accounting variable are ignored (equivalent to taking $\omega _{2}=..=\omega
_{n}=0).$ Our analysis identifies the locations of all the $n+1$ real roots
of $f(\kappa )=0$ in order to discuss the equation $f(\kappa )=\omega _{1}.$

Noting that $f(\lambda _{j}+)=-\infty $ and $f(\lambda _{j}-)=+\infty $, by
the Intermediate Value Theorem, there exist roots $\kappa _{j}(0)$ for $%
j=3,...,n$ of the equation $f(\kappa )=0$ which are real and satisfy%
\begin{equation*}
\lambda _{j}<\kappa _{j}(0)<\lambda _{j-1}.
\end{equation*}%
Finally, since $f(-\infty )=-\infty $ and $f(\lambda _{_{n}}-)=+\infty ,$
there exists a root $\kappa _{n+1}(0)<\lambda _{n}.$

Evidently the equation $f(\kappa )=\omega _{1}$ similarly has roots $\kappa
_{j}(\omega _{1})$ for any $\omega _{1}$ in the same intervals for $%
j=3,...,n+1$. Thus we have `interlacing' for $j\geq 2$:%
\begin{equation*}
\kappa _{n+1}(\omega _{1})<\lambda _{n}<\kappa _{n}(\omega _{1})<...<\lambda
_{j}<\kappa _{j}(\omega _{1})<\lambda _{j-1}<...<\kappa _{3}(\omega
_{1})<\lambda _{2}.
\end{equation*}

Now $f(\lambda _{2}+)=-\infty $ and $f(\lambda _{1})=0$ with $f(+\infty
)=-\infty .$ There are thus two generic possibilities.

(i) The equation $f(\kappa )=0$ has a root in $(\lambda _{2},\lambda _{1})$.
In this case, as $\omega _{1}$ increase from zero there are initially two
real roots $\kappa _{1}(\omega _{1})$ and $\kappa _{2}(\omega _{2})$ of the
equation $f(\kappa )=\omega _{1}$ which lie in $(\lambda _{2},\lambda _{1})$
and which move into coincidence. Thereafter they become complex conjugates
which behave qualitatively as in the case $n=2.$ In particular, $\kappa
_{3},...,\kappa _{n+1}$ are increasing in $\omega _{1}$ (under the
assumptions of this case), and since%
\begin{equation*}
\kappa _{1}+..+\kappa _{n+1}=\lambda _{1}+\lambda _{2}+...\lambda _{n}+\beta
,
\end{equation*}%
with the right-hand side constant, we see that $\Re(\kappa _{1})=\frac{%
1}{2}(\kappa _{1}+\kappa _{2})$ is decreasing in $\omega _{1}.$ As the lower
bound is $\frac{1}{2}(\lambda _{1}+\beta ),$ the complex roots certainly do
not enter the disc $|\zeta |\leq \lambda _{2}$ provided $\frac{1}{2}(\lambda
_{1}+\beta )\geq \lambda _{2}$ i.e. $\beta \geq 2\lambda _{2}-\lambda _{1}.$
(Note that $2\lambda _{2}-\lambda _{1}<\lambda _{1}.)$

(ii) The function $f(\kappa )$ is negative in $(\lambda _{2},\lambda _{1})$
and the equation $f(\kappa )=0$ has its remaining root in $(\lambda
_{1},+\infty ).$ In this case the convergence assumption is violated for
small $\omega _{1}$ in that there are eigenvalues of $H(\omega )$ outside
the disc $|\zeta |<\lambda _{1}.$

Finally, the special case arises when $f(\kappa )=0$ has a double root at $%
\kappa =\lambda _{1}.$ In this case as $\omega _{1}$ increases from zero the
remaining two roots are conjugate complex and again behave as in Figures $\ref{gig3a},\,\ref{gig3b}$.

\subsection{General $n:$ case $\protect\delta _{1}\protect\omega _{1}<0$
with $\protect\delta _{j}\protect\omega _{j}<0$ for $j=2,...,n.$}

This proceeds similarly. The root $\kappa _{1}$ decreases towards $\lambda
_{2}$ as $\omega _{1}$ increases. Likewise the roots $\kappa _{j}$ for $%
j=2,...,n-1$ decrease towards $\lambda _{j+1}.$ As these latter roots remain
bounded, the remaining two roots $\kappa _{n}$ and $\kappa _{n+1}$ may be
real, but will be the unbounded complex conjugates for large enough $\omega
_{1}$ (indeed the remaining component of the graph is $n$-shaped). This time
the real part increases towards $\frac{1}{2}(\beta +\lambda _{1})$ as $%
\omega _{1}$ increases.

\subsection{General $n:$ case when $\protect\delta _{k}\protect\omega _{k}=0$
for some $k=2,...,n$ with $\protect\delta _{j}\protect\omega _{j}>0$ for
remaining indices $j$}

We find that in all these cases DPI\ cannot hold at $R=\lambda _{2},$ by
Proposition 1.

Clearly, if $\delta _{2}=0,$ then $\kappa =\lambda _{2}$ is an eigenvalue of
the system, i.e. there is an eigenvalue in the annulus $\mathcal{A}$ and so
DPI\ cannot hold at $R=\lambda _{2}.$

In general, if $\delta _{k}\omega _{k}=0$ for just one of $k=3,...,n,$ then $%
\chi _{H}(\lambda _{k})=0,$ and in fact
\begin{eqnarray*}
\chi _{H}(\kappa ) &=&\prod\limits_{j=1}^{n+1}(\kappa -\lambda
_{j})-\sum_{j=1}^{n}\omega _{j}\delta _{j}\prod\limits_{h\neq j}^{n}(\kappa
-\lambda _{h}) \\
&=&(\kappa -\lambda _{k})[\prod\limits_{h=1,h\neq k}^{n+1}(\kappa -\lambda
_{h})-\sum_{j=1,j\neq k}^{n}\omega _{j}\delta _{j}\prod\limits_{h\neq
j}^{n}(\kappa -\lambda _{h})] \\
&=&(\kappa -\lambda _{k})\cdot \chi _{K}(\kappa ),
\end{eqnarray*}%
with $K$ the matrix obtained from $H$ by omitting the $k$-th row and column,
i.e. a reduction of order occurs. But now the assumptions for $K$ agree with
those made in the case of the previous subsection. It now follows that $\chi
_{K}(\kappa )$ has an eigenvalue in the annulus $\mathcal{A}$, and so again
DPI cannot hold at $R=\lambda _{2}.$

\bigskip

Finally, if $\delta _{k}\omega _{k}=0$ for several among $k=3,...,n$ and $%
\delta _{j}\omega _{j}>0$ for the remaining indices $j$, then a further
reduction of order occurs, with the same conclusion that an eigenvalue
exists in the annulus $\mathcal{A}$. So here too, DPI cannot hold at $%
R=\lambda _{2}.$

Note that in both scenarios the locus of $\kappa _{1}$ decreases as $\omega
_{1}$ increases from zero.

\subsection{Effect of the dividend-on-dividend multiplier: a distortion
example}

\begin{figure}[!t] %F1
\centering\includegraphics[width=0.7\textwidth] {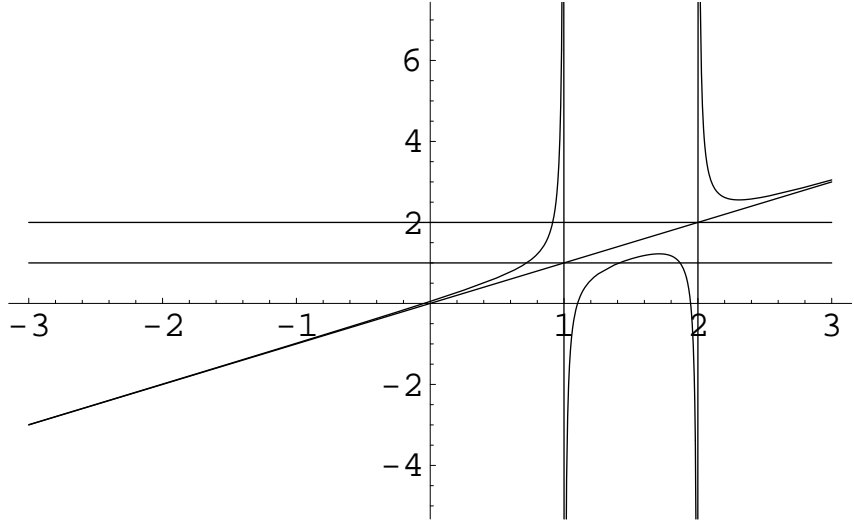}
\caption{
Graph of $\omega_3(\kappa)$ with $\delta_1=-1,\,\delta_2=+1,\,\omega_1,\,\omega_2>0$; $\kappa$-axis horizontal.
} \label{gig5a}
\end{figure}

We conclude this section by illustrating the effect of the policy parameter $%
\beta =\omega _{3}$ on the three eigenvalues in the case $\omega _{1}=\omega
_{2}=0.1.$ In the range $\beta <\lambda _{2}$ we see in the illustrative
example of Figure $\ref{gig5a}$ that the root $\kappa _{1}$ decreases whilst the root $%
\kappa _{2}$ increases as $\beta $ increases; $\kappa _{3}$ increases for
all $\beta ,$ as might be expected, with $\lambda _{2}$ as supremum.
Intuitively speaking, the push away from the origin created by the two
increasing roots $\kappa _{2}$ and $\kappa _{3}$ causes the location of the
coincident root $\kappa _{1}=\kappa _{2}$ to execute a jump up to a new
coincidence location above $\lambda _{1},$ by way of a continuous root locus
in the complex $\zeta $-plane (see the Remark on bifurcation in the next
section). The push can in fact be physically interpreted. The
partial-fraction-expansion terms in (\ref{beta}) may be regarded as
modelling electric charges placed at the pole locations $\lambda _{j}$ and
acting according to an inverse distance law (see Marden \cite[\S 1.3 p. 7]{Mar}).
Thus for $\beta $ large enough to ensure both $\kappa _{2}$ and $\kappa _{1}$
have been re-located above $\lambda _{1}$, we see the locus of $\kappa _{1}$
resume its downward path towards the origin (but tending in the limit only
as far as the barrier $\lambda _{1}$), while $\kappa _{2}$ resumes its
upward path away from the origin. The locus dynamics are investigated more
properly in the next section.%

\begin{figure}[!t] %F1
\centering\includegraphics[width=0.7\textwidth] {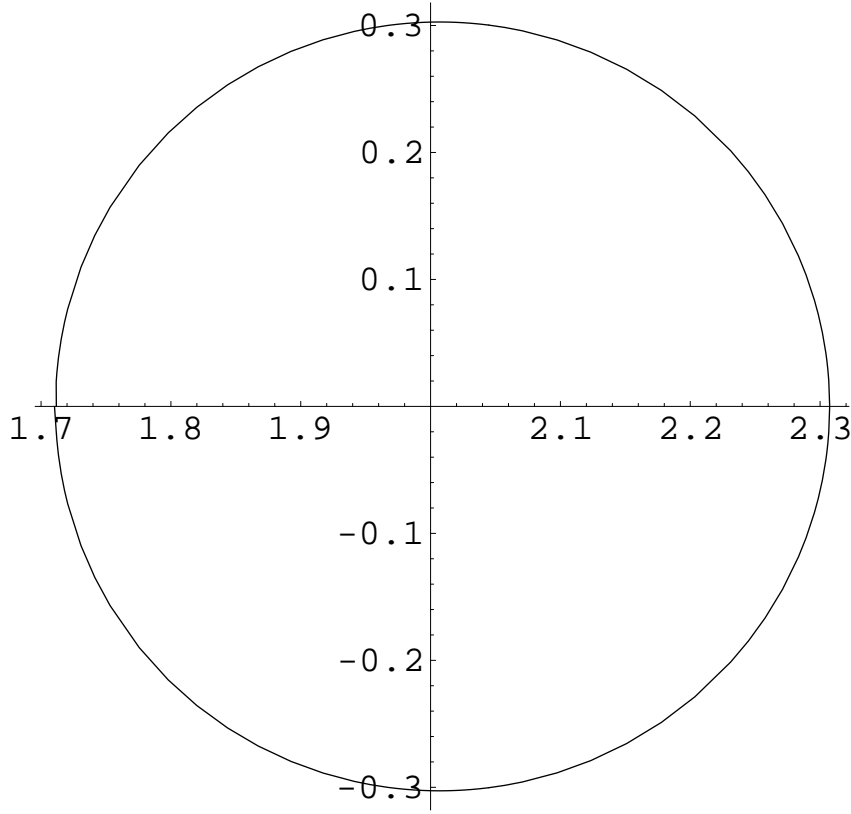}
\caption{
Root locus of conjugate roots as $\beta$ varies for $\omega_1,\omega_2>0$ fixed ( $\delta_1=-1,\,\delta_2=+1$).
} \label{gig5b}
\end{figure}

The conjugate root locus for this case is shown in Figure $\ref{gig5b}$. For $\beta $
varying over the range corresponding to the presence of complex roots (with $%
\omega _{1}$ and $\omega _{2}$ fixed), the root locus of the conjugate pair
appears close to being a circle $|\kappa -\lambda _{2}|^{2}=const.$ as a
consequence of the Circle Lemma and Proposition 6.

\subsubsection{A case study: distorted circle}

Since $\omega _{1}=\omega _{2}=0.1$ we have here
\begin{equation*}
f(\kappa )=10\kappa +\frac{1}{\kappa -2}-\frac{1}{\kappa -1},
\end{equation*}%
so the equation $f(x)=K$ is equivalent to $f_{2}(\kappa )=K-10\kappa .$ The
essential reason for the closeness observed here is the slow variation of $%
\kappa _{3}$ in the vicinity of $\lambda _{2}$. Indeed, by (\ref{poly}) we
have%
\begin{equation*}
\chi (\kappa )=(\kappa -\lambda _{1})(\kappa -\lambda _{2})(\kappa -\beta
)-\{\omega _{1}\delta _{1}(\kappa -\lambda _{2})+\omega _{2}\delta
_{2}(\kappa -\lambda _{1})\},
\end{equation*}%
from which we compute that%
\begin{eqnarray*}
\kappa _{1}+\kappa _{2}+\kappa _{3} &=&\lambda _{1}+\lambda _{2}+\beta , \\
\kappa _{1}\kappa _{2}\kappa _{3} &=&\lambda _{1}\lambda _{2}\beta -\omega
_{1}\delta _{1}\lambda _{2}-\omega _{2}\delta _{2}\lambda _{1},
\end{eqnarray*}%
and so
\begin{eqnarray*}
(\kappa _{1}-\lambda _{1})(\kappa _{2}-\lambda _{1}) &=&\kappa _{1}\kappa
_{2}-\lambda _{1}(\kappa _{1}+\kappa _{2})+\lambda _{1}^{2} \\
&=&(\lambda _{1}\lambda _{2}\beta -\omega _{1}\delta _{1}\lambda _{2}-\omega
_{2}\delta _{2}\lambda _{1})\kappa _{3}^{-1}-\lambda _{1}(\lambda _{1}+\beta
+(\lambda _{2}-\kappa _{3}))\\
&{\quad}&+\lambda _{1}^{2} \\
&=&\beta \lambda _{1}(\lambda _{2}-\kappa _{3})\kappa _{3}^{-1}-\lambda
_{1}(\lambda _{2}-\kappa _{3})-(\omega _{1}\delta _{1}\lambda _{2}+\omega
_{2}\delta _{2}\lambda _{1}\lambda _{2}^{-1}) \\
&=&\lambda _{1}(\lambda _{2}-\kappa _{3})[\beta \kappa _{3}^{-1}-1]-(\omega
_{1}\delta _{1}\lambda _{2}+\omega _{2}\delta _{2}\lambda _{1}\lambda
_{2}^{-1}).
\end{eqnarray*}%
So the positive sensitivity to variation in $\beta $ is slight, as asserted,
though in fact it is increasing.

Of course, writing $\lambda _{2}$ in place of $\kappa _{3},$ all the
dependence on $\beta $ is formally lost; however, note that the remaining
constant is inadequate as an explanataion of the actual radius.

\subsection{Qualitative behaviour in the most general case}

The picture emerging from our analysis for $\delta _{1}\omega _{1}<0$ is
that the dominant and subdominant eigenvalues remain in the annulus (of Th.
1) provided
\begin{equation*}
\beta >2\lambda _{1}-\lambda _{2}.
\end{equation*}
Other eigenvalues, when real, interlace between the eigenvalues of $A,$ with
$\kappa _{j}$ remaining associated with $\lambda _{h}$ for $h=j$ or $h=j+1.$
One pair of roots become unbounded as a complex conjugate pair and are
asymptotic on one or other side of $\Re(\zeta )=\frac{1}{2}(\lambda
_{1}+\beta )$ as $\omega _{1}\rightarrow +\infty ,$ depending on the
remaining parameters of $\omega .$ Evidently a regime change occurs with the
unbounded root locus being vertical. Other conjugate pairs execute deformed
circles as described by the Circle Lemma and Proposition 6 of \S 3.1.

\section{Differential properties of eigenvalues: some bifurcation analysis}

The purpose of this section is to analyse briefly the root locus in the $%
\zeta $-plane. We conduct a partial analysis mostly concentrated on the
dynamics of the dominant eigenvalue as $\omega _{1}$ changes (with the
remaining policy parameters fixed) with a view to completing the proof of
the dominance theorem in the next section. Our starting point is the
following proposition, which follows from (\ref{poly}) by implicit
differentiation. It is suggested by our earlier Circle Lemma of \S 3.1 and
Proposition 6 concerning the third derivative.

\bigskip

\noindent \textbf{Proposition 8} \textbf{\ (Under the assumption of distinct
eigenvalues). }\textit{In the canonical setting of Proposition 1, let the
dividend-policy vector be represented by }$\omega =(\omega _{1},...,\omega
_{n+1})$\textit{\ with }$\omega _{n+1}=\beta .$\textit{\ The eigenvalues }$%
\kappa _{j}=\lambda _{j}^{L}$\textit{\ of the augmented matrix }$%
A^{\lrcorner },$\textit{\ viewed} \textit{as functions of }$\omega =(\omega
_{1},...,\omega _{n+1}),$\textit{\ satisfy the following differential
properties expressed in terms of the constants }$\lambda _{j}=\lambda
_{j}^{A}$\textit{:}%

\medskip

\noindent\text{(i) for } $1\leq h\leq n+1,\text{ }1\leq k\leq n$
\begin{equation}
\frac{\partial \kappa _{h}}{\partial \omega _{k}}=-\frac{\delta _{j}(\lambda
_{1}-\kappa _{h})(\lambda _{2}-\kappa _{h})...(\lambda _{k-1}-\kappa
_{h})(\lambda _{k+1}-\kappa _{h})...(\lambda _{n}-\kappa _{h})}{(\kappa
_{1}-\kappa _{h})(\kappa _{2}-\kappa _{h})...(\kappa _{h-1}-\kappa
_{h})(\kappa _{h+1}-\kappa _{h})...(\kappa _{n+1}-\kappa _{h})};%
  \label{kiomegaj}
\end{equation}%
\noindent\text{(ii) for } $1\leq h\leq n+1$
\begin{equation}
\frac{\partial \kappa _{h}}{\partial \omega _{n+1}}=\frac{(\lambda
_{1}-\kappa _{h})(\lambda _{2}-\kappa _{h})...(\lambda _{n}-\kappa _{h})}{%
(\kappa _{1}-\kappa _{h})(\kappa _{2}-\kappa _{h})...(\kappa _{h-1}-\kappa
_{h})(\kappa _{h+1}-\kappa _{h})...(\kappa _{n+1}-\kappa _{h})}.
\label{kiomeganplus1}
\end{equation}

\bigskip

\begin{proof} From the identity%
\begin{equation*}
\chi _{H}(\kappa )=(\kappa -\kappa _{1} )(\kappa -\kappa _{2} )...(\kappa -\kappa
_{n+1} ),
\end{equation*}%
we have%
\begin{equation*}
\chi _{H}^{\prime }(\kappa )=\frac{d}{d\kappa }\chi _{H}(\kappa )=(\kappa-\kappa
_{2} )...(\kappa -\kappa _{n+1})+(\kappa -\kappa _{1})(\kappa -\kappa
_{3} )...(\kappa _{n+1}-\kappa )+...,
\end{equation*}%
so that%
\begin{equation*}
\chi _{H}^{\prime }(\kappa _{j})=\prod\limits_{h\neq j}(\kappa _{j}-\kappa
_{h}).
\end{equation*}%
So, for $\omega =\omega _{0}:=(0,...,0,\beta ),$ as then $\kappa_1=\lambda_1$,%
\begin{equation*}
\chi _{H}^{\prime }(\kappa _{1})=(\lambda _{1}-\lambda _{2})...(\lambda _{1}-\lambda
_{n})(\lambda _{1}-\omega _{n+1}).
\end{equation*}%
From the identity%
\begin{equation*}
\chi _{H}(\kappa _{j}(\omega ),\omega )=0,
\end{equation*}%
by implicit differentiation,%
\begin{equation*}
\frac{\partial \kappa _{1}}{\partial \omega _{j}}=-\frac{\partial \chi }{%
\partial \omega _{j}}\div \left( \frac{\partial \chi }{\partial \kappa }%
\right) _{\kappa =\kappa _{1}}=\frac{\partial \chi }{\partial \omega _{j}}%
\div \prod\limits_{h>1}(\kappa _{h}-\kappa _{1}),
\end{equation*}%
except at the critical points $\omega $ defined by%
\begin{equation*}
\prod\limits_{j>1}(\kappa _{j}(\omega )-\kappa _{1}(\omega ))=0,
\end{equation*}%
e.g. where the locus of $\kappa _{1}(\omega )$ crosses $\kappa _{2}(\omega
). $ The result of the Proposition now follows directly from (\ref{poly}). \qed
\end{proof}

\bigskip

\noindent \textit{Remark.} For the assumption of distinct roots to hold we
must manifestly disregard the non-generic critical points, which are those
points $\omega $ where any two of the functions $\kappa _{j}$ agree in
value; of particular importance to us are points $\omega $ where $\kappa
_{1}(\omega )$ may cease to be the largest eigenvalue (in modulus), as for
instance when it agrees in value with $\kappa _{2}(\omega )$. The first
formula when $j=1$ is to be read as%
\begin{equation*}
\frac{\partial \kappa _{1}}{\partial \omega _{1}}=-\frac{\delta _{1}(\lambda
_{2}-\kappa _{1})...(\lambda _{n}-\kappa _{1})}{(\kappa _{2}-\kappa
_{1})...(\kappa _{n}-\kappa _{1})(\kappa _{n+1}-\kappa _{1})},
\end{equation*}%
and note that, at $\omega =\omega _{0}:=(0,...,0,\beta ),$ we have, by (\ref%
{initialize}), $\kappa _{n+1}=\beta $ and for $j=1,...,n:$
\begin{equation*}
\kappa _{j}=\lambda _{j},
\end{equation*}%
\begin{equation}
\frac{\partial \kappa _{1}}{\partial \omega _{1}}=\frac{\delta _{1}}{%
(\lambda _{1}-\beta )},  \label{starter}
\end{equation}%
and%
\begin{equation*}
\frac{\partial \kappa _{1}}{\partial \omega _{j}}=0,\text{ for }j>1.
\end{equation*}%
The equation (\ref{starter}) implies that the choice of a $\beta $ value
close to $\lambda _{1}$ will accelerate the growth rate $\kappa _{1}$ of the
leading canonical variable $Z^{1}$ relative to the first dividend-policy
coefficient.

\bigskip

\noindent \textbf{Technical point}. In the arguments that follow, it is
important to realize that when the roots $\kappa _{j}$ and $\kappa _{j+1}$
are complex conjugates, then for real $\kappa _{1}$ the following signature
property is satisfied:%
\begin{equation*}
\text{sign}[(\kappa _{j}-\kappa _{1})(\kappa _{j+1}-\kappa _{1})]=+1,
\end{equation*}%
just as when $\kappa _{j}$ and $\kappa _{j+1}$ were real and both below $%
\kappa _{1}$. (Since the quadratic has no real roots, it is positive
definite here.)

\bigskip

\noindent \textbf{Corollary 1 (Bifurcation behaviour near }$\kappa
_{1}=\kappa _{2}$\textbf{). }\textit{Assume the eigenvalues of }$A$ \textit{%
are real and distinct and that for }$j=3,...,n,$\textit{\ }$\kappa _{j}$%
\textit{\ is real and satisfies }$\lambda _{j+1}<\kappa _{j}<\lambda _{j}.$
\textit{At any bifurcation point, for small enough positive increments in }$%
\omega _{1},$\textit{\ the conjugate complex roots }$\kappa _{1}$ \textit{%
and }$\kappa _{2}$\textit{\ move away from the origin if }$\delta _{1}>0$%
\textit{, and towards the origin if }$\delta _{1}<0$\textit{. }

\bigskip

\noindent \textit{Remark.} The corollary is thus in keeping with the
intuition expressed in connection with Figures $\ref{gig5a},\ref{gig5b}$ where we alluded to the
push away from the origin for the root $\kappa _{2}$ regarded as a function
of $\beta =\omega _{3}.$

\bigskip

\begin{proof} Suppose that\textbf{\ }$\delta _{1}=-1.$
Suppose that $\kappa _{1}=\kappa _{2}$ occurs at some point $\omega
_{1}=\omega _{1}^{\ast }.$ If now $\omega _{1}=\omega _{1}^{\ast }+\Delta
\omega $ with $\Delta \omega >0$, put $\kappa =\Re(\kappa _{1})$ and $%
\beta =\Im(\kappa _{1}),$ so that%
\begin{eqnarray*}
\kappa _{1} &=&\kappa +i\varepsilon ,\qquad \kappa _{2}=\kappa -i\varepsilon
, \\
(\kappa _{1}-\lambda _{j}) &=&\rho _{j}e^{i\theta _{j}},\qquad (\kappa
_{1}-\kappa _{j+1})=\rho _{j}^{\prime }e^{i\phi _{j}}.
\end{eqnarray*}%
We need to be sure which of $\kappa \pm i\varepsilon $ is to be interpreted
as $\kappa _{1}$ and $\kappa _{2}$ and whether the display above is correct.
In fact, either interpretation is valid, and leads to the same conclusion;
we return to this issue in a moment. Thus, since $\rho _{j}<\rho
_{j}^{\prime },$ for $\Delta \omega $ small enough we shall have%
\begin{equation*}
\theta _{j}>\phi _{j},
\end{equation*}%
so that%
\begin{equation*}
\frac{\rho _{j}e^{i\theta _{j}}}{\rho _{j}^{\prime }e^{i\phi _{j}}}=\frac{%
\rho _{j}}{\rho _{j}^{\prime }}\exp [i(\theta _{j}-\phi _{j})],
\end{equation*}%
and hence that%
\begin{equation*}
\frac{\Delta \kappa _{1}}{\Delta \omega }=-\frac{1}{2\varepsilon i}\frac{%
\rho _{2}e^{i\theta _{2}}}{\rho _{2}^{\prime }e^{i\phi _{2}}}\cdot ...\cdot
\frac{\rho _{n}e^{i\theta _{n}}}{\rho _{n}^{\prime }e^{i\phi _{n}}}=\frac{1}{%
2\varepsilon }\frac{\rho _{2}}{\rho _{2}^{\prime }}\cdot ...\cdot \frac{\rho
_{n}}{\rho _{n}^{\prime }}\exp [i(\frac{\pi }{2}+\psi )],
\end{equation*}%
where $\psi $ is small and positive. That is, the remaining ratios pull $%
\Delta \kappa _{1}$ in the same direction, towards the origin. Note that if
we switch the interpretation of $\kappa \pm i\varepsilon $ around, then the
angles $\theta _{j},\phi _{j}$ change sign, making $\psi $ small and
negative. However, the sign of $(\kappa _{2}-\kappa _{1})$ also switches.

In conclusion, the conjugate complex roots $\kappa _{1}$ and $\kappa _{2}$
initially move closer to the origin if $\delta _{1}<0.$ See Figures $\ref{gig5b},\ref{gig6}$ above
for an illustration in the case $n=2,$ where the third root of $\chi
_{A^{\lrcorner }}$ is evidently real. Note the vertical asymptote in the
complex $\zeta $-plane for $\Re(\zeta )=\frac{1}{2}(\lambda _{2}+\beta
)$ identified by Proposition 4.
\qed
\end{proof}

\bigskip

\noindent \textit{Remark (Bifurcation behaviour elsewhere).}\textbf{\ }%
Assuming that the first repeated root is not the dominant root, one may
attempt to repeat the argument at the other locations to observe a tug of
war between those ratios below the coincidence location pulling one way and
those above it pulling the other way. (We have noted in \S 4.5 the electric
force field interpretation.) Who wins this tug of war is determined by the
geometric considerations, and so we discover that there will be a critical
point $\lambda ,$ a watershed, such that to the right of $\lambda $ the
complex roots move towards the origin, whereas to the left they move away
from it.

\bigskip

\noindent \textbf{Corollary 2.} \textit{Assume all the eigenvalues }$\lambda
_{j}$ \textit{are real and positive, at some }$\omega $\textit{\ }$\kappa
_{1}$ \textit{is the maximal eigenvalue (in modulus) and is real, and further%
}%
\begin{equation*}
\lambda _{3}<\kappa _{1}<\lambda _{1}.
\end{equation*}%
\textit{Then}

\noindent (i) \ \textit{for }$j=1,2$%
\begin{eqnarray*}
\text{sign}[\partial \kappa _{1}/\partial \omega _{1}] &=&\text{sign}[\delta
_{1}]\text{sign}[(\kappa _{1}-\lambda _{2})], \\
\text{sign}[\partial \kappa _{1}/\partial \omega _{2}] &=&-\text{sign}%
[\delta _{2}],
\end{eqnarray*}

\noindent (ii) \ \textit{for }$j=3,...,n$%
\begin{equation*}
\text{sign}[\partial \kappa _{1}/\partial \omega _{j}]=-\text{sign}[\delta
_{j}]\text{sign}[(\kappa _{1}-\lambda _{2})],
\end{equation*}%
\textit{and finally for }$j=n+1$%
\begin{equation*}
\text{sign}[\partial \kappa _{1}/\partial \omega _{n+1}]=-1.
\end{equation*}

\bigskip

\begin{proof} Counting the signs gives
\begin{eqnarray*}
\text{sign}[(\kappa _{2}-\kappa _{1})(\kappa _{3}-\kappa _{1})...(\kappa
_{n+1}-\kappa _{1})]&=&(-1)^{n}\\
\text{sign}[(\lambda _{1}-\kappa _{1})(\lambda _{3}-\kappa _{1})...(\lambda
_{n}-\kappa _{1})]&=&(-1)^{n-2}.
\end{eqnarray*}
For $j\geq 3:$%
\begin{eqnarray*}
&\text{sign}&[(\lambda _{1}-\kappa _{1})(\lambda _{2}-\kappa _{1})...(\lambda
_{j-1}-\kappa _{1})(\lambda _{j+1}-\kappa _{1})...(\lambda _{n}-\kappa
_{1})]\\
&=&(-1)^{(n-3)}\text{sign}[(\lambda _{2}-\kappa _{1})] =(-1)^{(n-4)}\text{sign}[(\kappa _{1}-\lambda _{2})].
\end{eqnarray*}\qed
\end{proof}

\bigskip

\noindent \textbf{Corollary 3.} \textit{Suppose all the eigenvalues of }$A$
\textit{are real and positive, and }$\omega $ \textit{is small enough so that%
}%
\begin{equation*}
|\kappa _{n+1}|<...<|\kappa _{2}|<\kappa _{1},
\end{equation*}%
\textit{and }$\kappa _{j}$ \textit{is real with }%
\begin{equation*}
\lambda _{j+1}<\kappa _{j}<\lambda _{j-1}.
\end{equation*}%
\textit{Then, for each }$h$ \textit{with }$\kappa _{h}$\textit{\ real,}%
\begin{eqnarray*}
\text{sign}[\partial \kappa _{h}/\partial \omega _{h}] &=&\text{sign}[\delta
_{h}], \\
\text{sign}[\partial \kappa _{h}/\partial \omega _{j}] &=&\text{sign}[\delta
_{h}]\text{sign}[(\kappa _{h}-\lambda _{h})]\text{sign}[(\kappa _{h}-\lambda
_{j})],\text{\textit{\ }for }h\neq j\leq n, \\
\text{sign}[\partial \kappa _{h}/\partial \omega _{n+1}] &=&\text{sign}%
[(\kappa _{h}-\lambda _{h})],
\end{eqnarray*}%
\textit{where }$h>1.$ \textit{In particular, provided }$\kappa _{2}<\kappa
_{1},$\textit{\ }%
\begin{eqnarray*}
\text{sign}[\partial \kappa _{2}/\partial \omega _{1}] &=&-\text{sign}%
[\delta _{1}]\text{sign}[(\kappa _{2}-\lambda _{2})]\text{ ,\qquad sign}%
[\partial \kappa _{2}/\partial \omega _{2}]=\text{sign}[\delta _{2}], \\
\text{sign}[\partial \kappa _{2}/\partial \omega _{j}] &=&\text{sign}[\delta
_{j}]\text{sign}[(\kappa _{2}-\lambda _{2})]\text{ for }3\leq j\leq n, \\
\text{sign}[\partial \kappa _{2}/\partial \omega _{n+1}] &=&\text{sign}%
[(\kappa _{2}-\lambda _{2})].
\end{eqnarray*}

\bigskip

\begin{proof} Recalling for $h\leq n$ that%
\begin{equation*}
\frac{\partial \kappa _{h}}{\partial \omega _{h}}=-\frac{\delta _{h}(\lambda
_{1}-\kappa _{h})(\lambda _{2}-\kappa _{h})...(\lambda _{h-1}-\kappa
_{h})(\lambda _{h+1}-\kappa _{h})...(\lambda _{n}-\kappa _{h})}{(\kappa
_{1}-\kappa _{h})(\kappa _{2}-\kappa _{h})...(\kappa _{h-1}-\kappa
_{h})(\kappa _{h+1}-\kappa _{h})...(\kappa _{n+1}-\kappa _{h})},
\end{equation*}%
we compute that%
\begin{eqnarray*}
\text{sign}[(\kappa _{1}-\kappa _{h})(\kappa _{2}-\kappa _{h})...(\kappa
_{h-1}-\kappa _{h})(\kappa _{h+1}-\kappa _{h})...(\kappa _{n+1}-\kappa
_{h})] &=&(-1)^{(n+1-h)}, \\
\text{sign}[(\lambda _{1}-\kappa _{h})(\lambda _{2}-\kappa _{h})...(\lambda
_{h-1}-\kappa _{h})(\lambda _{h+1}-\kappa _{h})...(\lambda _{n}-\kappa
_{h})] &=&(-1)^{(n-h)},
\end{eqnarray*}%
and so%
\begin{eqnarray*}
&=&(-1)^{(n-h)}\text{sign}[(\lambda _{h}-\kappa _{h})] \\
&=&\text{sign}[(\lambda _{1}-\kappa _{h})(\lambda _{2}-\kappa
_{h})...(\lambda _{h-1}-\kappa _{h})(\lambda _{h}-\kappa _{h})(\lambda
_{h+1}-\kappa _{h})...(\lambda _{n}-\kappa _{h})] \\
&=&\text{sign}[(\lambda _{j}-\kappa _{h})]\text{sign}[(\lambda _{1}-\kappa
_{h})(\lambda _{2}-\kappa _{h})...(\lambda _{j-1}-\kappa _{h})(\lambda
_{j+1}-\kappa _{h})...(\lambda _{n}-\kappa _{h})].
\end{eqnarray*}\qed
\end{proof}

\section{Proof of the dominance theorem}

We may now put together the analysis of the last sections to deduce our main
result concerning the location of the dominant eigenvalue of $H(\omega )$.

\bigskip

\textbf{Proof of Eigenvalue Dominance Theorem.} We consider the first part
of the theorem only, as the more general result follows by a restatement of
the same argument. By selecting the sign of $\delta _{1}$ as $(-1)$ and of $%
\delta _{h}$ for $h>1$ as $(+1),$ we can arrange, given (\ref{starter}), for
the eigenvalue function $\kappa _{1}(\omega )$ identified by the condition $%
\kappa _{1}(\omega _{0})=\lambda _{1}$ to be decreasing in $\omega $ in the
region%
\begin{equation*}
\{\omega :\omega _{1}>0,...,\omega _{n+1}>0\},
\end{equation*}%
and so to remain below $\lambda _{1}.$ We have, however, to ensure that $%
\kappa _{1}(\omega )$ remains the maximal root. Recall that $\omega
_{0}=(0,...,0,\beta ).$ Since the remaining eigenvalues $\lambda _{j}=\kappa
_{j}(\omega _{0})$ are below $\lambda _{1}$ we may, by continuity, ensure
that the eigenvalues functions $\kappa _{2}(\omega ),...,\kappa
_{n+1}(\omega )$ of $H(\omega )$ also lie strictly below $\lambda _{1}$ and
that moreover $\kappa _{2}(\omega )<\kappa _{1}(\omega ).$ $\square $

\bigskip

\noindent \textit{Remark.}\textbf{\ }By Corollary 3 it is possible that,
following a path in parameter space, the locus of $\kappa _{2}$ intersects
that of $\kappa _{1}$. Note, however, that if upon intersection at $\omega
^{\ast }$ we were thereafter to have $\kappa _{1}(\omega ^{\prime })<\kappa
_{2}(\omega ^{\prime })$ for $\omega ^{\prime }$ close to $\omega ^{\ast },$
then provided the remaining eigenvalues remain below $\kappa _{1}(\omega
^{\ast }),$ the signs of all the derivatives $\partial \kappa _{1}/\partial
\omega _{j}$ and those of $\partial \kappa _{2}/\partial \omega _{j}$ would
switch, i.e. both loci would turn around, a contradiction. Thus, subject to
the assumption about the remaining eigenvalues, this implies that in fact $%
\omega ^{\ast }$ is at the boundary of that region in policy-parameter space
where $\kappa _{1}$ and $\kappa _{2}$ are both real. Moreover, according to (%
\ref{kiomegaj}) the graph has infinite slope at $\omega ^{\ast }.$ We
illustrate this point in the following simple example with $n=1,$ $\lambda
_{1}=1$ and $|\beta |<1.$

\bigskip

\begin{figure}[!t] %F1
\centering\includegraphics[width=0.7\textwidth] {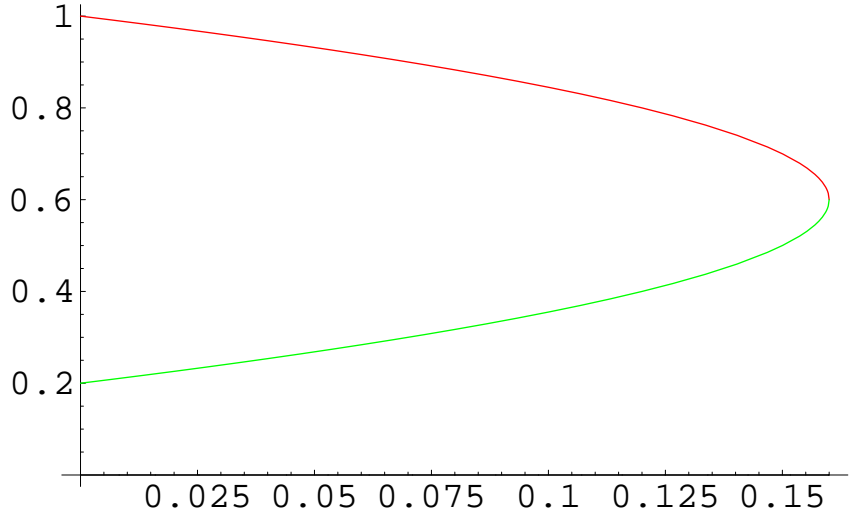}
\caption{
For fixed $\beta$ the root locus of $\kappa_1(\omega_1)$ (in red) and $\kappa_2(\omega_1)$ (in green) with vertical $\kappa$-axis.
} \label{gig6}
\end{figure}

\noindent \textbf{Example. }For $\delta _{1}=-1,$ let%
\begin{equation*}
A^{\lrcorner }=\left[
\begin{array}{cc}
1 & -1 \\
\omega _{1} & \beta%
\end{array}%
\right] .
\end{equation*}%
The characteristic polynomial is $\kappa ^{2}-(1+\beta )\kappa +(\beta
+\omega _{1}).$ Here
\begin{equation*}
\frac{d\kappa _{1}}{d\omega }=\frac{1}{(1+\beta -2\kappa _{1})}=\frac{-b_{1}%
}{(\kappa _{2}-\kappa _{1})},
\end{equation*}%
since $\kappa _{1}+\kappa _{2}=1+\beta .$ The roots are real for $$\omega
_{1}\leq \frac{1}{4}[(1+\beta )^{2}-4\beta ]=\frac{1}{4}(1-\beta
)^{2}=\omega _{1}^{\ast }(\beta ),$$ and we have
\begin{equation*}
\kappa _{1}=\frac{1}{2}\left( 1+\beta +\sqrt{(1-\beta )^{2}-4\omega _{1}}%
\right)
\end{equation*}

decreasing down to $\kappa _{1}=\frac{1}{2}(1+\beta )$ as $\omega _{1}$
increases, and analogously
$$\kappa _{2}=\frac{1}{2}\left( 1+\beta -\sqrt{(1-\beta )^{2}-4\omega _{1}}\right)$$
increasing up to $\kappa _{2}=\frac{1}{%
2}(1+\beta ).$ See Figure $\ref{gig6}$.

We note that as the roots become complex the real part stays constant at $%
(1+\beta )/2,$ i.e. the root locus bifurcates and the conjugate roots move
orthogonally to the real axis; there being no poles in this simple case,
there is no `push' on the real part, neither away nor towards the origin. In
particular, provided the roots are in the unit circle, they remain in the
annulus $\beta =\lambda _{2}<|\zeta |<1=\lambda _{1}.$

\section{Obtaining dividend irrelevance}

In this section we prove the results in Proposition 1.

After a change of accounting state variables from $z_{t}$ to, say $Z_{t},$
the system becomes%
\begin{equation*}
\left[
\begin{array}{c}
Z_{t+1} \\
d_{t+1}%
\end{array}%
\right] =H\left[
\begin{array}{c}
Z_{t} \\
d_{t}%
\end{array}%
\right] ,
\end{equation*}%
where

\begin{equation*}
H=H(\omega )=\left[
\begin{array}{ccccc}
\lambda _{1} & 0 &  & 0 & b_{1} \\
0 & \lambda _{2} &  &  & b_{2} \\
&  & ... &  &  \\
0 & 0 &  & \lambda _{n} & b_{n} \\
\omega _{1} & \omega _{2} &  &  & \omega _{n+1}%
\end{array}%
\right] ,
\end{equation*}%
with $\omega _{n+1}=\beta ,$ the new augmented matrix (with the same
eigenvalues as the original augmented matrix $A^{\lrcorner })$ and where $%
\lambda _{1},...,\lambda _{n}$ are the eigenvectors of $A$ assumed presented
in decreasing modulus size (with $\lambda _{1}$ largest).

Evidently the characteristic polynomial
\begin{equation*}
\chi _{H}(\kappa )=\chi _{H}(\kappa ,\omega
_{1},...,\omega _{n},\omega _{n+1})=|\kappa I-H|
\end{equation*}%
is the same as $\chi _{A^{\lrcorner }}(\kappa )$. We assume its eigenvalues $%
\kappa _{1},...,\kappa _{n+1}$ have distinct modulus.

The force of our implicit assumptions on the \textbf{dividend significance
coefficients} is that they are all non-zero: $b_{j}\neq 0$ for all $j.$
Otherwise, the eigenvalues of the augmented matrix would include all the
eigenvalues of the reduced matrix. (To see this expand the characteristic
determinant by the $j$-th row.)

Henceforth we assume the canonical variables have been re-scaled by $b_{j}$
and we may therefore take for the \textbf{canonical dividend significance
coefficients} the symbol $\delta _{j}$ with the additional stipulation that%
\begin{equation*}
|\delta _{j}|=1\text{ for }j=1,..,n.
\end{equation*}

As a first step we note the consequence for dividend irrelevance of the
non-zero dividend significance coefficients. Writing $Z=(...,Z^{j},...)$ and
fixing $j,$ there are coefficients $l_{1},..,l_{n+1}$ such that%
\begin{equation*}
d_{t}=\sum\nolimits_{h}l_{h}\kappa _{h}^{t}.
\end{equation*}%
Now the equation%
\begin{equation*}
Z^{j}(t+1)=\lambda _{j}Z^{j}(t)+\delta _{j}\sum\nolimits_{h}l_{h}\kappa
_{h}^{t},
\end{equation*}%
with the solution also given by the eigenvalues of the augmented matrix%
\begin{equation*}
Z^{j}(t)=\sum\nolimits_{h}L_{h}\kappa _{h}^{t},
\end{equation*}%
must satisfy%
\begin{equation*}
\sum L_{h}\kappa _{h}^{t}(\kappa _{h}-\lambda _{j})=\delta
_{j}\sum\nolimits_{h}l_{h}\kappa _{h}^{t}
\end{equation*}%
for all $t.$ Hence%
\begin{equation*}
L_{h}=\frac{\delta _{j}l_{h}}{\kappa _{h}-\lambda _{j}}.
\end{equation*}%
We thus have, assuming $R>|\kappa _{j}|$ for all $j,$ that the dividend
series converges, and
\begin{equation*}
P_{0}(R;d)=\sum_{t=1}^{\infty }R^{-t}d_{t}=\sum_{h}\frac{1}{R-\kappa _{h}}%
l_{h}\kappa _{h}=\frac{1}{\delta _{j}}\sum_{h}\frac{\kappa _{h}-\lambda _{j}%
}{R-\kappa _{h}}L_{h}\kappa _{h}.
\end{equation*}%
Consequently, if $R=\lambda _{j}$ is permitted, then
\begin{equation*}
P_{0}(R;d)=-\frac{1}{\delta _{j}}\sum_{h}L_{h}\kappa _{h}=-\frac{Z_{1}^{j}}{%
\delta _{1}}=-\frac{\lambda _{j}Z_{0}^{j}+\delta _{j}d_{0}}{\delta _{1}}=-%
\frac{\lambda _{j}Z_{0}^{j}}{\delta _{1}}+d_{0}.
\end{equation*}%
This indeed depends only on the initial data. See Ashton [Ash] for a
discussion of this formula. The earliest form of this equation is due to
Ohlson in 1989, though published later in \cite{Ohl96}.

We recall that the basis of this calculation is the identity%
\begin{eqnarray*}
P_{0}(R;d) &=&\sum_{t=1}^{\infty }R^{-t}d_{t}=\sum_{t=1}^{\infty
}R^{-t}\sum_{h}l_{h}\kappa _{h}^{t}=\sum_{h}l_{h}\sum_{t=1}^{\infty
}R^{-t}\kappa _{h}^{t} \\
&=&\sum_{h}l_{h}\frac{\kappa _{h}/R}{1-\kappa _{h}/R}=\sum_{h}\frac{1}{%
R-\kappa _{h}}l_{h}\kappa _{h}.
\end{eqnarray*}

We will thus obtain dividend irrelevance at $R=\lambda _{j}$ provided all
the eigenvalues $\kappa _{j}$ are in modulus less than $R.$

\section{Derivation of equivalences}

We begin by expanding by the bottom row%
\begin{equation*}
|H-\kappa I|=\left\vert
\begin{array}{ccccc}
\lambda _{1}-\kappa & 0 &  & 0 & \delta _{1} \\
0 & \lambda _{2}-\kappa &  &  & \delta _{2} \\
&  & ... &  &  \\
0 & 0 &  & \lambda _{n}-\kappa & \delta _{n} \\
\omega _{1} & \omega _{2} &  &  & \lambda _{n+1}-\kappa%
\end{array}%
\right\vert _{n+1}=0
\end{equation*}%
to obtain
\begin{equation*}
(-1)^{n}\omega _{1}D_{1}(\kappa )-...-\omega _{n}D_{n}(\kappa
)+\prod\limits_{j=1}^{n+1}(\lambda _{j}-\kappa )=0,
\end{equation*}%
or%
\begin{equation*}
(-1)^{n}\omega _{1}D_{1}(\kappa )+(-1)^{n-1}\omega _{2}D_{2}(\kappa
)...-\omega _{n}D_{n}(\kappa
)=(-1)(-1)^{n+1}\prod\limits_{j=1}^{n+1}(\kappa -\lambda _{j}),
\end{equation*}%
where%
\begin{eqnarray*}
D_{1}(\kappa ) &=&\left\vert
\begin{array}{cccc}
0 & 0 & 0 & \delta _{1} \\
\lambda _{2}-\kappa &  & 0 & \delta _{2} \\
... &  &  &  \\
0 &  & \lambda _{n}-\kappa & \delta _{n}%
\end{array}%
\right\vert _{n} \\
&=&(-1)^{n-1}\delta _{1}\prod\limits_{j=2}^{n}(\lambda _{j}-\kappa )=\delta
_{1}\prod\limits_{j=2}^{n}(\kappa -\lambda _{j}).
\end{eqnarray*}

Similarly,%
\begin{eqnarray*}
D_{2}(\kappa ) &=&\left\vert
\begin{array}{ccccc}
\lambda _{1}-\kappa & 0 &  & 0 & \delta _{1} \\
0 & 0 &  & 0 & \delta _{2} \\
& \lambda _{3}-\kappa &  &  &  \\
&  &  &  &  \\
0 & 0 &  & \lambda _{n}-\kappa & \delta _{n}%
\end{array}%
\right\vert _{n} \\
&=&(-1)^{n-2}\delta _{2}\prod\limits_{j\neq 2}^{n}(\lambda _{j}-\kappa
)=-\delta _{2}\prod\limits_{j\neq 2}^{n}(\kappa -\lambda _{j}).
\end{eqnarray*}%
This yields the equation%
\begin{equation*}
\omega _{1}\delta _{1}\prod\limits_{j=2}^{n}(\kappa -\lambda _{j})+\omega
_{2}\delta _{2}\prod\limits_{j\neq 2}^{n}(\kappa -\lambda _{j})+...+\omega
_{n}\delta _{n}\prod\limits_{j\neq n}^{n}(\kappa -\lambda
_{j})=\prod\limits_{j=1}^{n+1}(\kappa -\lambda _{j}).
\end{equation*}%
Dividing by $\prod\limits_{h\neq j}^{n}(\kappa -\lambda _{h})$,
\begin{equation*}
(\kappa -\lambda _{j})(\kappa -\beta )=\omega _{j}\delta _{j}+\sum_{h\neq
j}\omega _{h}\delta _{h}\frac{\kappa -\lambda _{j}}{\kappa -\lambda _{h}}%
=\omega _{j}\delta _{j}+\sum_{h\neq j}\omega _{h}\delta _{h}\left( 1-\frac{%
\lambda _{j}-\lambda _{h}}{\kappa -\lambda _{h}}\right) ,
\end{equation*}%
or%
\begin{equation*}
\omega _{1}\delta _{1}+\omega _{2}\delta _{2}+...=(\kappa -\lambda
_{j})(\kappa -\beta )+\sum_{h\neq j}^{n}\frac{\omega _{h}\delta _{h}(\lambda
_{j}-\lambda _{h})}{\kappa -\lambda _{h}},
\end{equation*}%
as required.

Dividing by $\prod\limits_{j=1}^{n}(\kappa -\lambda _{j}),$ we obtain%
\begin{equation*}
\kappa -\beta =\frac{\omega _{1}\delta _{1}}{\kappa -\lambda _{1}}+\frac{%
\omega _{2}\delta _{2}}{\kappa -\lambda _{2}}+...+\frac{\omega _{n}\delta
_{n}}{\kappa -\lambda _{n}}.
\end{equation*}

\section{Invertible parametrization and zero placement}

This section is devoted to a proof of Proposition 3. Let us write%
\begin{equation*}
\chi _{A}(\kappa )=|\kappa I-A|=\sum_{s=0}^{n}(-1)^{s}a_{s}\kappa
^{n-s}=\prod\limits_{j=1}^{n}(\kappa -\lambda _{j}),
\end{equation*}%
so that%
\begin{equation*}
a_{s}=\sum_{j_{1}<...<j_{s}}\lambda _{j_{1}}...\lambda _{j_{s}}
\end{equation*}%
denotes the elementary symmetric function summing the zeros of $\chi
_{A}(\kappa )$ taken $s$ at a time. Thus%
\begin{equation*}
a_{0}=1,\qquad a_{1}=\lambda _{1}+...+\lambda _{n},\qquad ...\qquad
a_{n}=\lambda _{1}...\lambda _{n}.
\end{equation*}%
Hence%
\begin{eqnarray*}
\prod\limits_{j=1}^{n+1}(\kappa -\lambda _{j}) &=&(\kappa -\beta )\left[
\sum_{s=0}^{n}(-1)^{s}a_{s}\kappa ^{n-s}\right] \\
&=&\kappa ^{n+1}-(a_{1}+\beta a_{0})\kappa ^{n}+...+(-1)^{s}[a_{s+1}+\beta
a_{s}]\kappa ^{n-s}\\
&&+...+(-1)^{n+1}\beta a_{n}.
\end{eqnarray*}%
As a first step we compute that%
\begin{equation*}
\lambda _{2}+...+\lambda _{n}=a_{1}-\lambda _{1},
\end{equation*}%
and that%
\begin{eqnarray*}
\sum_{1<u<v}\lambda _{u}\lambda _{v} &=&\sum_{u<v}\lambda _{u}\lambda
_{v}-\lambda _{1}\sum_{1<v}\lambda _{v}=\sum_{u<v}\lambda _{u}\lambda
_{v}-\lambda _{1}(\sum_{v}\lambda _{v}-\lambda _{1}) \\
&=&a_{1}-\lambda _{1}(a_{1}-\lambda )=a_{2}-a_{1}\lambda _{1}+\lambda
_{1}^{2}.
\end{eqnarray*}%
Similarly,%
\begin{eqnarray*}
\sum_{1<u<v<w}\lambda _{u}\lambda _{v}\lambda _{w} &=&\sum_{u<v<w}\lambda
_{u}\lambda _{v}\lambda _{w}-\lambda _{1}\sum_{1<v<w}\lambda _{v}\lambda _{w}
\\
&=&a_{3}-\lambda _{1}[a_{2}-a_{1}\lambda _{1}+\lambda _{1}^{2}] \\
&=&a_{3}-a_{2}\lambda _{1}+a_{1}\lambda _{1}^{2}-\lambda _{1}^{3}.
\end{eqnarray*}%
The pattern is now clear, and we shall show by induction that%
\begin{equation*}
\sum_{1<j_{1}<...<j_{s}}\lambda _{j_{1}}...\lambda
_{j_{s}}=a_{s}-a_{s-1}\lambda _{1}+\lambda
_{1}^{2}a_{s-2}+...+(-1)^{s}a_{0}\lambda _{1}^{s}.
\end{equation*}%
Indeed,%
\begin{eqnarray*}
\sum_{j_{1}<...<j_{s}}\lambda _{j_{1}}...\lambda _{j_{s}}
&=&\sum_{j_{1}<...<j_{s}}\lambda _{j_{1}}...\lambda _{j_{s}}-\lambda
_{1}\sum_{1<j_{2}<...<j_{s}}\lambda _{j_{2}}...\lambda _{j_{s}} \\
&=&a_{s}-\lambda _{1}(a_{s-1}+...+(-1)^{s-1}\lambda _{1}^{s-1}) \\
&=&a_{s}-a_{s-1}\lambda _{1}+\lambda
_{1}^{2}a_{s-2}+...+(-1)^{s}a_{0}\lambda _{1}^{s}.
\end{eqnarray*}

Note that%
\begin{equation*}
a_{n}-a_{n-1}\lambda _{1}+\lambda _{1}^{2}a_{n-2}+...+(-1)^{n}a_{0}\lambda
_{1}^{n}=0,
\end{equation*}%
so%
\begin{equation*}
\lambda _{2}...\lambda _{n}=\lambda _{1}^{-1}a_{n}=a_{n-1}-\lambda
_{1}a_{n-2}+...+(-1)^{n-1}\lambda _{1}^{n-1}.
\end{equation*}%
Our next step is to observe that the coefficients in the polynomial on the
right-hand side of identity (\ref{poly}) may be expanded as follows:%
\begin{eqnarray*}
D(\kappa ) &=&\omega _{1}\delta _{1}\prod\limits_{j\neq 1}^{n}(\kappa
-\lambda _{j})+\omega _{2}\delta _{2}\prod\limits_{j\neq 2}^{n}(\kappa
-\lambda _{j})+...+\omega _{n}\delta _{n}\prod\limits_{j\neq n}^{n}(\kappa
-\lambda _{j}) \\
&=&(\omega _{1}\delta _{1}+\omega _{2}\delta _{2}+...+\omega _{n}\delta
_{n})\kappa ^{n-1}-(\omega _{1}\delta _{1}[\lambda _{2}+...]+...)\kappa
^{n-2} \\
&&+(\omega _{1}\delta _{1}[\lambda _{2}\lambda _{3}+...]+...)\kappa ^{n-3}+
\\
&&...+(-1)^{s}(\omega _{1}\delta _{1}\bar{a}_{s}(1)+...)\kappa
^{n-s}+...+(-1)^{n-1}[\sum_{j=1}^{n}\omega _{j}\delta
_{j}\prod\limits_{h\neq j}^{n}\lambda _{h}],
\end{eqnarray*}%
where%
\begin{equation*}
\bar{a}_{s}(h)=\sum_{\substack{ j_{1}<...<j_{s}  \\ j_{k}\neq h}}\lambda
_{j_{1}}...\lambda _{j_{s}}
\end{equation*}%
(i.e. the summation refers to the omission of $h$ from any of the components
$j_{1}...j_{s}$). Note also that%
\begin{equation*}
\prod\limits_{h\neq j}^{n}\lambda _{h}=\frac{a_{n}}{\lambda _{j}}.
\end{equation*}

We now consider for constants $p_{s}$ the identity
\begin{align*}
\prod\limits_{j=1}^{n+1}(\kappa -\lambda _{j})-\{\omega _{1}\delta
_{1}\prod\limits_{j\neq 1}^{n}(\kappa -\lambda _{j})+\omega _{2}\delta
_{2}\prod\limits_{j\neq 2}^{n}(\kappa -\lambda _{j})+...+\omega _{n}\delta
_{n}\prod\limits_{j\neq n}^{n}(\kappa -\lambda _{j})\} \\
=\kappa ^{n+1}-p_{0}\kappa ^{n}+p_{2}\kappa ^{n-1}+...+p_{n}.
\end{align*}

Comparing sides, we obtain%
\begin{eqnarray*}
p_{0} &=&(a_{1}+\beta a_{0}), \\
p_{1} &=&a_{2}+\beta a_{1}-(\omega _{1}\delta _{1}+\omega _{2}\delta
_{2}+...+\omega _{n}\delta _{n}), \\
p_{2} &=&a_{3}+\beta a_{2}-\left( \omega _{1}\delta _{1}\lambda
_{1}+...\right) , \\
&&... \\
p_{n} &=&a_{n+1}+\beta a_{n}-\left( \omega _{1}\delta _{1}\lambda
_{1}^{-1}+...\right) .
\end{eqnarray*}%
Now given any $p_{0},$ we select $\beta ,$ so that%
\begin{equation*}
\beta =p_{0}-a_{1}.
\end{equation*}%
For the remaining equations, we have%
\begin{eqnarray*}
p_{1}-a_{2}-\beta a_{1} &=&\omega _{1}\delta _{1}+\omega _{2}\delta
_{2}+...+\omega _{n}\delta _{n}, \\
p_{2}-a_{3}-\beta a_{2} &=&\omega _{1}\delta _{1}(a_{1}-\lambda _{1})+..., \\
p_{3}-a_{4}-\beta a_{3} &=&\omega _{1}\delta _{1}(a_{2}-s_{1}\lambda
_{1}+\lambda _{1}^{2})+..., \\
&&... \\
p_{n}-a_{n+1}-\beta a_{n} &=&\omega _{1}\delta _{1}(a_{n-1}-a_{n-2}\lambda
_{1}+(-1)^{n-1}a_{0}\lambda _{1}^{n-1})+...,
\end{eqnarray*}%
where $a_{n+1}=0,$ i.e.%
\begin{equation*}
N\left[
\begin{array}{c}
\delta _{1}\omega _{1} \\
\delta _{2}\omega _{2} \\
... \\
\delta _{n}\omega _{n}%
\end{array}%
\right] =\left[
\begin{array}{c}
p_{1}-a_{2}-\beta a_{1} \\
p_{2}-a_{3}-\beta a_{2} \\
... \\
p_{n}-a_{n+1}-\beta a_{n}%
\end{array}%
\right] .
\end{equation*}%
Here the coefficient matrix $N$ is given as follows:%
\begin{equation*}
N=\left[
\begin{array}{ccc}
1 &  & 1 \\
a_{1}-\lambda _{1} &  & a_{1}-\lambda _{n} \\
a_{2}-a_{1}\lambda _{1}+\lambda _{1}^{2} & ... & a_{2}-a_{1}\lambda
_{n}+\lambda _{n}^{2} \\
... &  & ... \\
a_{n-1}-\lambda _{1}a_{n-2}+...+(-1)^{n-1}\lambda _{1}^{n-1} &  &
a_{n-1}-\lambda _{n}a_{n-2}+...+(-1)^{n-1}\lambda _{n}^{n-1}%
\end{array}%
\right] .
\end{equation*}%
Its determinant is equal, up to a possible sign change, to the van der Monde
determinant $V(\lambda _{1},...,\lambda _{n}).$ Hence $N$ is non-singular
and the equation may be solved for any given vector $(p_{1},...,p_{n}).$ To
see this, note that $N$ may be reduced to the alternant matrix $A(0,...,n-1)$
in the variables $(-\lambda _{1}),...(-\lambda _{n}):$
\begin{eqnarray*}
&&\left\vert
\begin{array}{cccc}
1 & 1 & ... & 1 \\
a_{1}-\lambda _{1} & a_{1}-\lambda _{2} & ... & a_{1}-\lambda _{n} \\
a_{2}-a_{1}\lambda _{1}+\lambda _{1}^{2} &  &  &  \\
... &  &  &  \\
\sum_{i=0}^{n-1}(-1)^ia_{n-1-i}\;\lambda_{1}^i &  &  &
\sum_{i=0}^{n-1}(-1)^ia_{n-1-i}\;\lambda_{n}^i
\end{array}%
\right\vert \\[2mm]
&&\qquad\qquad =\left\vert
\begin{array}{cccc}
1 & 1 & ... & 1 \\
-\lambda _{1} & -\lambda _{2} & ... & -\lambda _{n} \\
\lambda _{1}^{2} &  &  & \lambda _{n}^{2} \\
... &  &  &  \\
(-\lambda _{1})^{n-1} &  &  & (-\lambda _{n})^{n-1}%
\end{array}%
\right\vert \\[2mm]
&&\qquad \qquad=V(-\lambda _{1},...,-\lambda _{n})
\end{eqnarray*}%
(taking $a_{1}$ times the first row, $a_{2}$ times the second row and so on).

It is now easy to find the inverse transformation by applying the elementary
row operations just used to the original matrix equation. This leads to the
following result. Putting $g_{j}=p_{j}-a_{j+1}-\beta a_{j},$ the original
equations
\begin{equation*}
N\left[
\begin{array}{c}
\delta _{1}\omega _{1} \\
\delta _{2}\omega _{2} \\
... \\
\delta _{n}\omega _{n}%
\end{array}%
\right] =\left[
\begin{array}{c}
g_{1} \\
g_{2} \\
... \\
g_{n}%
\end{array}%
\right]
\end{equation*}%
now transform to%
\begin{equation*}
V\left[
\begin{array}{c}
\delta _{1}\omega _{1} \\
\delta _{2}\omega _{2} \\
... \\
\delta _{n}\omega _{n}%
\end{array}%
\right] =\left[
\begin{array}{c}
h_{1} \\
h_{2} \\
\\
\\
h_{n}%
\end{array}%
\right] ,
\end{equation*}%
where%
\begin{eqnarray*}
h_{1} &=&g_{1}, \\
h_{2} &=&g_{2}-a_{1}h_{1}, \\
h_{3} &=&g_{3}-a_{2}h_{1}+a_{1}h_{2}, \\
&&... \\
h_{n} &=&g_{n}-a_{n-1}h_{1}+a_{n-2}h_{2}-...\pm a_{1}h_{n-1}.
\end{eqnarray*}%
Note that with the sign adjustment $h_{n}^{\prime }=(-1)^{n}h_{n}$ the
equations specify $g_{j}$ as a convolution (taking $a_{0}:=1).$ Now we have%
\begin{equation*}
\left[
\begin{array}{c}
\delta _{1}\omega _{1} \\
\delta _{2}\omega _{2} \\
... \\
\delta _{n}\omega _{n}%
\end{array}%
\right] =V^{-1}\left[
\begin{array}{c}
h_{1} \\
h_{2} \\
\\
\\
h_{n}%
\end{array}%
\right] ,
\end{equation*}%
where the inverse $V^{-1}$ is given by (see Klinger \cite{Klin}) the matrix with $%
jk$ entry%
\begin{equation*}
(-1)^{j+k}\;\frac{\bar{a}_{n-j}(k)}{\prod\limits_{l=1}^{k-1}(\lambda
_{k}-\lambda _{l})\prod\limits_{h=k+1}^{n}(\lambda _{h}-\lambda _{k})};
\end{equation*}%
here $\bar{a}_{s}(j)$ is the elementary symmetric function as above (sum
over the $s$-fold product omitting the variable $\lambda _{j}$).

\section{Asymptotics of the unbounded roots}

\textbf{Proof of Proposition 4}. Assume that $\kappa \neq \lambda _{j}$ for $%
j=1,...,n.$ We rewrite (\ref{omega1}) as%
\begin{eqnarray*}
\omega _{1}\delta _{1}+\omega _{2}\delta _{2}+...+\omega _{n}\delta _{n}
&=&(\kappa -\lambda _{j})(\kappa -\beta )+\sum_{k\neq j}^{n}\frac{\omega
_{k}\delta _{k}(\lambda _{j}-\lambda _{k})}{\kappa } \\
&&+\sum_{k\neq j}^{n}\omega _{k}\delta _{k}(\lambda _{j}-\lambda _{k})\left[
\frac{1}{\kappa -\lambda _{k}}-\frac{1}{\kappa }\right] ,
\end{eqnarray*}%
so that%
\begin{eqnarray*}
\omega _{1}\delta _{1}+\omega _{2}\delta _{2}+...+\omega _{n}\delta _{n}
&=&(\kappa -\lambda _{j})(\kappa -\beta )+\frac{1}{\kappa }\left[
\sum_{k\neq j}^{n}\omega _{k}\delta _{k}(\lambda _{j}-\lambda _{k})\right] \\
&&+\sum_{k\neq j}^{n}\omega _{k}\delta _{k}(\lambda _{j}-\lambda _{k})\left[
\frac{\lambda _{\kappa }}{\kappa (\kappa -\lambda _{k})}\right] .
\end{eqnarray*}%
Iterating, we obtain%
\begin{eqnarray}\label{identity-N}
\omega _{1}\delta _{1}+\omega _{2}\delta _{2}+...+\omega _{n}\delta _{n}
&=&(\kappa -\lambda _{j})(\kappa -\beta )\\
&&+\sum_{s=1}^{N}\frac{1}{\kappa ^{s}}%
\left[ \sum_{k\neq j}^{n}\omega _{k}\delta _{k}\lambda _{\kappa
}^{s}(\lambda _{j}-\lambda _{k})\right]  \notag\\
&&+\sum_{k\neq j}^{n}\omega _{k}\delta _{k}\lambda _{\kappa }^{N-1}(\lambda
_{j}-\lambda _{k})\frac{1}{\kappa ^{N}(\kappa -\lambda _{k})}.  \notag
\end{eqnarray}%
and hence the assertion of the Proposition follows provided $\kappa
>|\lambda _{k}|.$ The alternative direct derivation by expanding $\kappa
^{-1}(1-\lambda _{k}/\kappa )^{-1}$as a geometric series is less informative
about the convergence of the series.

We may in principle use the identity (\ref{identity-N}) (valid for all large
enough $\kappa $) recursively to obtain an asymptotic expansion (in $\omega
_{j})$ for the unbounded roots.

To obtain the first term of the expansion, let $|\kappa |\rightarrow \infty $
and consider the quadratic approximation%
\begin{eqnarray*}
\omega _{1}\delta _{1}+\omega _{2}\delta _{2}+...+\omega _{n}\delta _{n}
&\sim &\kappa ^{2}-(\lambda _{j}+\beta )\kappa +\lambda _{j}\beta \\
&=&[\kappa -\frac{1}{2}(\lambda _{j}+\beta )]^{2}-\frac{1}{4}(\lambda
_{j}-\beta )^{2},
\end{eqnarray*}%
so that $\omega _{j}\delta _{j}$ is large and positive and
\begin{equation*}
\kappa =\frac{1}{2}(\lambda _{j}+\beta )\pm \sqrt{\frac{1}{4}(\lambda
_{j}-\beta )^{2}+\sum_{k}\omega _{k}\delta _{k}}.
\end{equation*}%
Hence for $\delta _{j}=-1$ we have the first term of the expansion to be
\begin{equation*}
\kappa =\frac{1}{2}(\lambda _{j}+\beta )\pm i\sqrt{|\omega _{j}|}.
\end{equation*}

Now write
\begin{equation*}
\kappa =\alpha +i\sqrt{|\omega _{j}|},\qquad \alpha =\hat{\kappa}%
_{j}+\varepsilon ,\qquad \hat{\kappa}_{j}=\frac{1}{2}(\lambda _{j}+\beta
),\qquad A=\sum_{k\neq j}^{n}\omega _{k}\delta _{k}(\lambda _{j}-\lambda
_{k}).
\end{equation*}%
With this notation, taking only one term in (\ref{identity-N}) we have with $%
\omega _{j}$ large and negative that%
\begin{eqnarray*}
-\omega _{j}+\sum_{k\neq j}^{n}\omega _{k}\delta _{k} &=&[\kappa -\hat{\kappa%
}_{j}]^{2}-\frac{1}{4}(\lambda _{j}-\beta )^{2}+\frac{A}{\kappa } \\
&=&[i\sqrt{\omega _{j}}+\varepsilon ]^{2}-\frac{1}{4}(\lambda _{j}-\beta
)^{2}+\frac{A}{\alpha +i\sqrt{\omega _{j}}} \\
&=&-\omega _{j}+2i\varepsilon \sqrt{\omega _{j}}+\varepsilon ^{2}-\frac{1}{4}%
(\lambda _{j}-\beta )^{2}-A\frac{\alpha -i\sqrt{\omega _{j}}}{\omega
_{j}-\alpha ^{2}},
\end{eqnarray*}%
so to first order in $\varepsilon $%
\begin{equation*}
\sum_{k\neq j}^{n}\omega _{k}\delta _{k}=2i\varepsilon \sqrt{\omega _{j}}-%
\frac{1}{4}(\lambda _{j}-\beta )^{2}-A\frac{\hat{\kappa}_{j}+\varepsilon -i%
\sqrt{\omega _{j}}}{\omega _{j}-\hat{\kappa}_{j}^{2}}.
\end{equation*}%
Writing $\varepsilon =u+iv$ and taking real and imaginary parts gives the
two equations%
\begin{eqnarray*}
\frac{-A}{\omega _{j}-\hat{\kappa}_{j}^{2}}u-2v\sqrt{\omega _{j}} &=&\gamma
_{j}=\sum_{k\neq j}^{n}\omega _{k}\delta _{k}+\frac{1}{4}(\lambda _{j}-\beta
)^{2}+A\frac{\hat{\kappa}_{j}}{\omega _{j}-\hat{\kappa}_{j}^{2}}, \\
2u\sqrt{\omega _{j}}-\frac{A}{\omega _{j}-\hat{\kappa}_{j}^{2}}v &=&-\frac{A%
\sqrt{\omega _{j}}}{\omega _{j}-\hat{\kappa}_{j}^{2}}.
\end{eqnarray*}%
Note that for $|\omega _{j}|\rightarrow \infty $%
\begin{equation*}
\gamma _{j}\rightarrow \sum_{k\neq j}^{n}\omega _{k}\delta _{k}+\frac{1}{4}%
(\lambda _{j}-\beta )^{2}.
\end{equation*}

Now the determinant of the two equations in $u$ and $v$ above is positive
and equal to%
\begin{equation*}
\Delta =\frac{A^{2}}{(\omega _{j}-\hat{\kappa}_{j}^{2})^{2}}+4\omega _{j}.
\end{equation*}%
Solving for $u$ and $v,$ we have%
\begin{eqnarray*}
u\Delta &=&-\frac{A\gamma _{j}}{\omega _{j}-\hat{\kappa}_{j}^{2}}-\frac{%
2A\omega _{j}}{\omega _{j}-\hat{\kappa}_{j}^{2}}=-2A-\frac{2\hat{\kappa}%
_{j}^{2}+\gamma _{j}}{\omega _{j}-\hat{\kappa}_{j}^{2}}A, \\
v\Delta &=&\frac{-A\sqrt{\omega _{j}}}{(\omega _{j}-\hat{\kappa}_{j}^{2})^{2}%
}-2\gamma _{j}\sqrt{\omega _{j}},
\end{eqnarray*}%
so that%
\begin{equation*}
u=-\frac{A}{2\omega _{j}}+O(\omega _{j}^{-2}),\qquad v=-\frac{\gamma _{j}}{4%
\sqrt{\omega _{j}}}+O(\omega ^{-5/2}).
\end{equation*}%
The result for $u$ assumes that $A\neq 0$ and further that%
\begin{equation*}
2\hat{\kappa}_{j}^{2}+\gamma _{j}\neq 0,
\end{equation*}%
i.e.%
\begin{equation*}
\sum_{k\neq j}^{n}\omega _{k}\delta _{k}+\frac{3}{4}(\lambda _{j}+\beta
)^{2}\neq 0.
\end{equation*}

\section{Strip-and-two-circles theorem}

We prove Proposition 5. For part (i) we argue as follows. Suppose, for all $%
j $, that $\omega _{j}\delta _{j}\geq 0$ and $\beta \leq \lambda _{n}.$
Suppose $z$ satisfies
\begin{equation*}
z-\beta =\frac{\omega _{1}\delta _{1}}{z-\lambda _{1}}+\frac{\omega
_{2}\delta _{2}}{z-\lambda _{2}}+...+\frac{\omega _{n}\delta _{n}}{z-\lambda
_{n}}.
\end{equation*}%
If $z$ is strictly to the right of $\lambda _{1},$ then we may also assume
that $z$ has positive imaginary part (otherwise switch to the conjugate root
$\bar{z}$). The argument of $z-\lambda _{j}$ is thus positive for each $j,$
and that of $1/(z-\lambda _{j})$ negative, i.e. has negative imaginary part.
The right-hand side therefore sums to a complex number with negative
imaginary part. However, $z-\beta $ has positive imaginary part.

If $z$ is to the left of $\beta ,$ then we may suppose it has negative
imaginary part. The argument of $\lambda _{j}-z$ is for each $j$ thus
positive, as also for $\beta -z.$ Now apply the previous reasoning to the
identity%
\begin{equation*}
\beta -z=\frac{\omega _{1}\delta _{1}}{\lambda _{1}-z}+\frac{\omega
_{2}\delta _{2}}{\lambda _{2}-z}+...+\frac{\omega _{n}\delta _{n}}{\lambda
_{n}-z}.
\end{equation*}

For part (ii), let $\varepsilon _{j}$ be arbitrary real for $j=1,...n.$ We
will apply Marden's `Mean-Value Theorem for polynomials' (Marden \cite[\S %
2.8 p.23]{Mar}) to the polynomials $h_{j}$ for $j=1,...,n$ and the polynomial $%
f(z)$ as defined by
\begin{equation*}
f(z)=\prod\limits_{h=1}^{n+1}(z-\lambda _{h}),\qquad h_{j}(z)=\varepsilon
_{j}\prod\limits_{h\neq j}^{n}(z-\lambda _{h}).
\end{equation*}%
We must, however, first find for each $j$ the location of the roots of the
equation $f(z)=h_{j}(z).$ The roots are of course $z=\lambda _{k}$ for $%
k\neq j$ taken together with the two real roots of%
\begin{equation*}
(z-\beta )(z-\lambda _{j})=\varepsilon _{j},
\end{equation*}%
which are to the left of $\beta $ and the right of $\lambda _{j}$. The exact
and approximate formulas are%
\begin{equation*}
u_{j}^{\pm }=\frac{(\beta +\lambda _{j})\pm \sqrt{(\lambda _{j}-\beta
)^{2}+4\varepsilon _{j}}}{2}\sim \beta -\frac{\varepsilon _{j}}{4(\lambda
_{j}-\beta )},\qquad \lambda _{j}+\frac{\varepsilon _{j}}{4(\lambda
_{j}-\beta )},
\end{equation*}%
and require that%
\begin{equation*}
-\frac{1}{4}(\lambda _{j}-\beta )^{2}\leq \varepsilon _{j}.
\end{equation*}%
Thus the roots of all the equations lie in the interval $%
K=(u_{1}^{-},u_{1}^{+}).$ By Marden's Theorem in the special case of real
positive scalars $m_{j}$ summing to unity, the roots of%
\begin{equation*}
f(z)=\sum m_{j}h_{j}(z)
\end{equation*}%
lie in the star-shaped region $S(K,\pi /(n+1))$ (cf. Prop. 5). Thus if we
take $\varepsilon _{j}=\varepsilon $ small and $m_{j}\varepsilon =\delta
_{j}\omega _{j}$ so that
\begin{equation*}
\delta _{1}\omega _{1}+...+\delta _{n}\omega _{n}=\varepsilon ,\text{ with }%
\delta _{1}\omega _{1},...,\delta _{n}\omega _{n}\geq 0,
\end{equation*}%
then indeed $\sum m_{j}=1$ and all the roots of (\ref{char}) lie in the said
star-shaped region.

In fact, one may take $\varepsilon _{1}=\varepsilon $ small and $%
m_{1}\varepsilon _{1}=\delta _{1}\omega _{1},$ and for $j>1,$ $\varepsilon
_{j}=\mu :=\min \{\lambda _{j}-\lambda _{j+1}:$ for $j>1\}$ and $m_{j}\mu
=\delta _{j}\omega _{j}>0,$ leading to the restriction%
\begin{equation*}
1=m_{1}+...+m_{n}=\frac{\delta _{1}\omega _{1}}{\varepsilon }+\frac{1}{\mu }%
(\delta _{2}\omega _{2}+...),
\end{equation*}%
i.e.
\begin{equation*}
\delta _{1}\omega _{1}+\frac{\varepsilon }{\mu }(\delta _{2}\omega
_{2}+...+\delta _{n}\omega _{n})=\varepsilon \qquad (\delta _{1}\omega
_{1},...,\delta _{n}\omega _{n}\geq 0).
\end{equation*}

\section{The third-derivative test}

The result in \S 3.3 is a consequence of the following by specializing $%
g(\kappa ):=f(\kappa ^{\ast }+\kappa )-f(\kappa ^{\ast })$ when $f^{\prime
}(\kappa ^{\ast })=0.$

\bigskip

\noindent \textbf{Proposition 9.} \textit{For }$g$ \textit{with }$%
g(0)=0,g^{\prime }(0)=0,$\textit{\ }$g^{\prime \prime }(0)>0$\textit{, }$%
g^{\prime \prime \prime }(0)\neq 0,$ \textit{the solution }$\kappa =\kappa
(\omega )$ \textit{over the complex domain of the equation}%
\begin{equation*}
g(\kappa )=-\omega ^{2}
\end{equation*}%
\textit{with }$\omega >0$\textit{\ small\ and subject to }$\kappa (0)=0,$
\textit{satisfies }%
\begin{eqnarray*}
\Re(\kappa (\omega ))\mathit{\ initially\ increasing\ if\ }g^{\prime
\prime \prime }(0) &>&0,\mathit{\ and} \\
\mathit{initially\ decreasing\ if\ }g^{\prime \prime \prime }(0) &<&0.
\end{eqnarray*}%
\textit{The imaginary part is initially increasing and satisfies }$|\Im%
(\kappa (\omega ))|>\omega .$

\bigskip

\begin{proof} Below we work to order $o(\omega )$. Without loss
of generality we assume that $g^{\prime \prime }(0)=+2.$ (Otherwise rescale $%
g$ and $\omega ^{2}$ by $2/g^{\prime \prime }(0).$) Set $G:=g^{\prime \prime
\prime }(0)/6\neq 0;$ then, by Taylor's Theorem, the first approximation to
the equation is%
\begin{equation*}
g(\kappa )=\kappa ^{2}+o(\kappa ^{2})=-\omega ^{2},
\end{equation*}%
with solution $\kappa =\pm i\omega $. So we introduce correction terms by
putting
\begin{equation*}
\kappa =\alpha +i(\omega +\beta )
\end{equation*}%
and solve a second approximation%
\begin{equation*}
g(\kappa )=\kappa ^{2}+G\kappa ^{3}+o(\kappa ^{3})=-\omega ^{2}.
\end{equation*}%
Substitution for $\kappa $ gives after cancellation of the term $-\omega
^{2} $%
\begin{align} \label{re-and-im}
0=[\alpha ^{2}-(2\omega \beta +\beta ^{2})+2i\alpha (\omega +\beta
)]+G[\alpha ^{3}+3i\alpha ^{2}(\omega +\beta )\\
-3\alpha (\omega +\beta
)^{2}-i(\omega +\beta )^{3}]. \notag
\end{align}%
Equate real and imaginary parts to $0;$ cancelling the second by $(\omega
+\beta )$ (non-zero, w.l.o.g.) gives

\begin{equation*}
3G\alpha ^{2}+2\alpha -G(\omega +\beta )^{2}=0.
\end{equation*}%
The roots corresponding to the imaginary part, for $\omega $ and $\beta $
sufficiently small, are%
\begin{align*}
\frac{-1\pm \sqrt{1+3G^{2}(\omega +\beta )^{2}}}{3G}=&\;\frac{1}{3G}\{\frac{1}{2%
}3G^{2}(\omega +\beta )^{2}-\frac{1}{8}9G^{4}(\omega +\beta )^{4}+..\}\\
=&\;\frac{1}{2}G(\omega +\beta )^{2}+...,
\end{align*}%
So we may neglect $\alpha ^{2}$ and $\beta ^{2}$ in what follows, which
leads to the approximation%
\begin{equation*}
\alpha =G\omega \beta .
\end{equation*}\qed
\end{proof}

\bigskip

\textbf{Claim.} \textit{For }$\omega >0$\textit{\ small enough, }$\beta >0.$
\textit{In particular, }$\alpha $ \textit{is the same sign as }$G$\textit{.
(So, also, conversely, if }$\alpha $\textit{\ is the same sign as }$G,$%
\textit{\ then }$\beta >0.)$

\bigskip

\begin{proof} The equation for the real part of (\ref{re-and-im}),
ignoring $(2\omega \beta +\beta ^{2})$ and cancelling by $\alpha \neq 0,$
gives
\begin{equation*}
G\alpha ^{2}+\alpha -3G(\omega +\beta )^{2}=0.
\end{equation*}%
Its two roots $\alpha _{\pm }$ have negative product $-3(\omega +\beta
)^{2}, $ and%
\begin{eqnarray*}
\alpha _{\pm } &=&-\frac{1}{2G}\left[ 1\pm \sqrt{1+12G^{2}(\omega +\beta
)^{2}}\right] \\
&=&-\frac{1}{2G}\left[ 1\pm \{1+6G^{2}(\omega +\beta )^{2}\}%
\right] +o(\omega ) \\
&=&-G^{-1},\qquad 3G(\omega +\beta )^{2}\text{ to order }o(\omega ).
\end{eqnarray*}%
The root near $-1/G$ is ruled out by the continuity of the root locus, which
tends to $0$ as $\omega \rightarrow 0.$ For $G>0$ the positive root is near $%
3G(\omega +\beta )^{2}\sim 6G\omega \beta ,$ and so $\beta >0,$ as $\omega
>0.$ For $G<0,$ the negative root is $3G(\omega +\beta )^{2}\sim 6G\omega
\beta ,$ and again $\beta >0,$ as $\omega >0.$ \qed
\end{proof}

\bigskip

\noindent Mathematics Department, London School of Economics, Houghton Street, London
WC2A 2AE; A.J.Ostaszewski@lse.ac.uk
\end{document}